\documentclass[11pt]{amsart} 
\usepackage{a4wide}
\usepackage{amsmath,amsthm, amssymb}
\usepackage{fancybox}
\usepackage{mathbbol}
\usepackage{times}
\usepackage[alphabetic]{amsrefs}
\usepackage{hyperref}

\theoremstyle{plain}
\newtheorem*{main theorem}{Main Theorem}
\newtheorem*{theorem*}{Theorem}
\newtheorem*{emptytheorem*}{}
\newtheorem{theorem}{Theorem}
\newtheorem{lemma}{Lemma}[section]

\newtheorem*{proposition*}{Proposition}

\newtheorem{conjecture}{Conjecture}
\newtheorem*{conjecture*}{Conjecture}
\theoremstyle{definition}
\newtheorem{definition}{Definition}
\newtheorem*{definition*}{Definition}
\theoremstyle{remark}

\newtheorem*{remark*}{Remark}

\newcommand{\dist}{\text{dist}}

\makeindex
\begin{document}
\title[Nonuniformly expanding maps]{Stochastic-like behaviour \\ in  
nonuniformly expanding maps}

\author{Stefano Luzzatto}
\address{Mathematics Dept., Imperial College,  London}
\email{\href{mailto:stefano.luzzatto@imperial.ac.uk}{stefano.luzzatto@imperial.ac.uk}}

\urladdr{\href{http://www.ma.ic.ac.uk/~luzzatto}
{http://www.ma.ic.ac.uk/\textasciitilde luzzatto}}
\date{14 May 2004} 

\thanks{Thanks to Gerhard Keller, Giulio Pianigiani, Mark Pollicott, 
Omri Sarig, and 
Benoit Saussol for 
useful comments on a preliminary version of this paper.  Simon 
Cedervall made some useful contributions to some parts of section 
\ref{indMark}.}
\maketitle

 \setcounter{tocdepth}{2}
 \tableofcontents
\section{Introduction}

\subsection{Determinism versus randomness}
A feature of many real-life phenomena in areas as diverse as 
physics, biology, finance, economics, and many others, is the 
\emph{random-like} behaviour of processes which nevertheless 
are clearly \emph{deterministic}. On the level of applications 
this dual aspect has proved very problematic. Specific 
mathematical models tend to be developed either on the basis that the 
process is deterministic, in which case 
sophisticated numerical techniques can be used to attempt to 
understand and predict the evolution, or that it is random, in which 
case probability theory is used to model the process. 
Both approaches  lose sight of what is probably the 
most important and significant characteristic of the system which is 
precisely that 
it is deterministic \emph{and} has random-like behaviour. 
The theory of Dynamical Systems has 
contributed a phenomenal amount of work showing 
 that it is perfectly natural for completely deterministic 
systems to behave in a very random-like way and achieving a quite 
remarkable understanding of the mechanisms by which this occurs. 
The purpose of these 
notes is to survey some of this research. 

\subsection{Nonuniform expansivity}
\index{Nonuniform expansivity}
\index{Nonuniformly expanding}
We shall assume that the state space can be represented by 
a compact Riemannian manifold \( M \) and that the evolution of the 
process is given by a map 
\( f: M \to M \)  which is piecewise differentiable. 
Following an approach  which goes back at least to the 
first half of the 20th century, we shall discuss  how 
certain statistical properties can be deduced from 
geometrical assumptions on \( f \) formulated explicitly in terms of 
assumptions on
the \emph{derivative map} \( Df \) of \( f \). 
The basic strategy is to construct certain 
geometrical structures which then imply some 
statistical/probabilistic properties of the dynamics (a striking 
and pioneering example of this is the work of Hopf on the 
ergodicity of geodesic flows on manifolds of  negative 
curvature  \cite{Hop39}). 
Research work following this basic line of reasoning goes 
under the heading of \emph{Hyperbolic Dynamics} and/or \emph{Smooth 
Ergodic Theory}. 

The main focus of these notes will be on maps which satisfy an 
asymptotic expansivity condition. 

\begin{definition} 
    We say that \( f: M \to M \) is 
    \emph{(nonuniformly) expanding}
    if there exists \( \lambda>0 \) such that 
    \begin{equation*} \tag*{\( (*) \)}
    \liminf_{n\to\infty}\frac{1}{n}
    \sum_{i=0}^{n-1}\log \|Df^{-1}_{f^{i}(x)}\|^{-1}>\lambda
    \end{equation*}
    for almost every \( x\in M \). 
Equivalently, for 
almost every \( x \in M\) there exists a constant \( C_{x}>0 \) such 
that 
    \[ 
    \prod_{i=0}^{n-1} \|Df^{-1}_{f^{i}(x)}\|^{-1} \geq 
C_{x}e^{\lambda  n}
    \]
    for every \( n\geq 1 \).
\end{definition} 
The definition and the corresponding results can be 
generalized to the case in which the expansivity condition 
holds only on an invariant set of positive measure instead of on the 
entire manifold \( M \). See also \cite{Alv03} for a detailed 
treatment 
of the theory of nonuniformly expanding maps. 
Notice that the condition \( \log \|Df^{-1}_{x}\|^{-1}>0 \), which 
is equivalent to  \( \|Df^{-1}_{x}\|^{-1}>1 \) and which in turn 
is equivalent to \( \|Df^{-1}_{x}\|<1 \), implies that 
\emph{all} vectors in all directions are contracted by the
inverse of \( Df_{x} \) and thus that \emph{all} vectors in all 
directions are
\emph{expanded} by \( Df_{x} \); the intuitively more obvious 
condition \( \log \|Df_{x}\|>0 \), which is equivalent to 
\( \|Df_{x}\|>1 \), implies only that there is \emph{at
least} one direction in which vectors are expanded by \( Df_{x} \).
Thus a map is nonuniformly expanding if 
every vector is asymptotically expanded at a uniform 
exponential rate. The constant \( C_{x} \) can in 
principle be arbitrarily small and indicates that an arbitrarily 
large number of 
iterates may be needed before this exponential growth becomes 
apparent. 

In the special case in which condition \( (*) \) holds at 
\emph{every} point \( x \) and the constant \( C \) can be chosen 
uniformly positive independent of \( x \) we say that \( f \) is 
\emph{uniformly expanding}.  Thus \emph{uniformly expanding} is a 
special case of \emph{nonuniformly expanding}. The terminology is 
slightly awkward for historical reasons: uniformly expanding maps 
have 
traditionally been referred to simply as expanding maps whereas this 
term should more appropriately 
refer to the more general (i.e. possibly nonuniformly) expanding 
case. 
We shall generally say that \( f \) is strictly 
nonuniformly expanding if \( f \) satisfies condition \( (*) \) but 
is strictly \emph{not} uniformly expanding.
A basic theme of these notes is to 
discuss the difference between uniformly and nonuniformly expanding 
maps: how the nonuniformity affects the results and the ideas and 
techniques used in the proofs and how different \emph{degrees of 
nonuniformity} can me quantified. 

Nonuniform expansivity 
is  a special case 
of  \emph{nonuniform hyperbolicity}.   
\index{Nonuniform hyperbolicity}
This concept was first formulated and studied 
by Pesin \cites{Pes76, Pes77} and has since 
become one of the main areas of research in dynamical systems, see 
\cites{Bunetal89, You95h, BarPes02} and 
\cite{BarPes04} by Barreira and Pesin in this volume, 
for extensive and in-depth surveys. 
The formal definition is 
in terms of \emph{non-zero Lyapunov exponents} which 
means that the tangent 
bundle can be decomposed into subbundles in which vectors either 
contract or expand at an asymptotically exponential rate. 
Nonuniform \emph{expansivity } 
corresponds to the  
case in which all the Lyapunov exponents are positive and therefore 
all vectors expand asymptotically at an exponential rate. 
The natural setting for this situation is that 
 of (non-invertible) local diffeomorphisms  whereas the theory of 
nonuniform hyperbolicity has been developed mainly for 
diffeomorphisms (however
see also \cite{Rue82} and \cite{BarPes04}*{Section 5.8}).
For greater generality, and also because this has great importance for
applications, we shall also allow 
various kinds of critical and/or singular points for 
 \( f \) or its derivative. 

\subsection{General overview of the notes}
We first review the basic notions of invariant measure, 
ergodicity, mixing, and decay of correlations in order to fix the 
notation and to motivate the results and techniques. 
In section \ref{markov} we discuss the key idea of a \emph{Markov 
Structure} and sketch some of the arguments 
used to study systems which admit 
such a structure.  In section \ref{expanding}, \ref{almost expanding} 
and \ref{critical}, we give a 
historical and  technical survey of many classes of 
systems for which results are known, giving 
references to the original proofs whenever possible, and sketching in 
varying amounts of details the construction of Markov structures in 
such systems. In section \ref{general} we 
present some recent abstract results which go towards a general 
theory of 
nonuniformly expanding maps. In section \ref{existence} we discuss 
the important problem of verifying the geometric nonuniform 
expansivity assumptions in specific classes of maps. 
Finally, in 
section \ref{final} we make some concluding remarks and present some 
open 
questions and conjectures.

The focus on Markov structures is 
partly a matter of personal preference; 
in some cases  the results can be proved and/or were first proved 
using 
completely different arguments and techniques. 
Of particular importance is the so-called \emph{Functional-Analytic} 
approach in which the problems are re-formulated and 
reduced to questions about the spectrum of a 
certain linear operator on some functional space. There are several 
excellent survey texts focussing on this approach, see \cites{Bal01, 
Liv04, Via}. 
In any case, it is hard to see how the study of systems in which the 
hyperbolicity or expansivity is \emph{nonuniform}
can be carried out without constructing or defining some kind of 
subdivision  into subsets on which relevant 
estimates satisfy  uniform bounds. 
The Markov structures to be described below provide one 
very useful way in which this can be done and give 
some concrete geometrical structure. 
It seems very likely that these structures will prove 
useful in studying many other features of nonuniformly hyperbolic or 
expanding systems such as their stability, persistence, and even 
existence in particular settings. 
Another quite different way to partition a 
set satisfying nonuniform hyperbolicity conditions is with 
so-called \emph{Pesin} or \emph{regular} sets, see 
\cite{BarPes04}*{Section 4.5}. These sets play a very useful role in 
the general theory of nonuniform hyperbolicity for diffeomorphisms, 
for example in the construction of the stable and unstable 
foliations.

We shall always assume that \( M \) is a smooth, 
compact, 
Riemannian manifold of dimension \( d \geq 1 \). For 
simplicity we shall call the Riemannian volume \emph{Lebesgue 
measure}, denote it by \( m \) or \( |\cdot | \) and assume that it 
is normalized so 
that \( m(M) = |M|= 1 \).   We let \( f: M \to M \) denote a 
Lebesgue-measurable map. In practice we shall always assume 
significantly more regularity on \( f \), e.g. that \( f \) is 
\( C^{2} \) or at least piecewise \( C^{2} \), but the main 
definitions apply in the more general case of \( f \) measurable. 
 All measures on \( M \) will be assumed to be defined on the 
Borel \( \sigma \)-algebra of \( M \).

\subsection{Invariant measures}
\index{Invariant measure}
For a set \( A \in M \) and a map \( f: M \to M \) we define \( 
f^{-1}(A) = \{x: f(x) \in A\}. \)
\begin{definition}
We say that a probability measure \( \mu \) on \( M \) is 
\emph{invariant} under  
 \( f \) if 
\[ 
\mu(f^{-1}(A)) = \mu(A)
 \]
for every \( \mu \)-measurable set \( A \subset M \).
\end{definition}
A given measure can be invariant for many different maps. For example 
Lebesgue measure on the circle \( M = \mathbb S^{1} \) is invariant 
for 
the identity map \( f(\theta) = \theta \), the rotation \( f(\theta) 
= 
\theta + \alpha \) for any \( \alpha\in\mathbb R \), and the covering 
map \( f(\theta) = \kappa \theta \) for any \( \kappa \in\mathbb N 
\).  Similarly, for a given point \( p\in M  \), 
the Dirac-\( \delta \) measure \( \delta_{p} \) 
defined by 
\[
\delta_p(A) = \left\{ \begin{array}{ll} 
		   1 & \mbox{$p \in A$}\\
		   0 & \mbox{$p \notin A$}
		  \end{array} 
		  \right .
\]
is invariant for any map \( f \) for which \( f(p) = p \). On the 
other hand, a given map \( f \) can admit many invariant measures. 
For example \emph{any} probability measure is invariant for 
the identity map \( f(x) = x \) and, more generally, any map which 
admits multiple fixed or periodic points admits as invariant measures 
the Dirac-\( \delta \) measures supported on such fixed points or 
their natural generalizations distributed along the orbit of the 
periodic points. 
There exist also maps that do not admit any invariant 
probability measures. However some mild conditions, e.g. continuity 
of \( f \), do guarantee that there exists at least one. 

A first step in the application of the theory and methods of ergodic
theory is to introduce some ways of distinguishing 
between the various invariant measures.  We do this by 
introducing various properties which such measures may or may not 
satisfy. Unless we specify otherwise we shall use \( \mu \) to denote 
a generic invariant probability measure for a given unspecified map 
\( 
f: M \to M \). 

\subsection{Ergodicity}
\index{Ergodicity}
\begin{definition}
    We say that \( \mu \) is \emph{ergodic} 
if there does  not exist a measurable set \( A \) with 
 \[
 f^{-1}(A) = A   \quad\text{ and }\quad  \mu(A) \in (0,1) 
    \]
    \end{definition}
In other words, any fully invariant set \( A \), i.e. a set for which 
\( 
f^{-1}(A) = A \), has either zero or full measure. This is a kind of 
\emph{indecomposability} property of the measure. If such a set 
existed, its complement \( B=A^{c} \) would also be fully invariant 
and , in particular, both \( A \) and \( B \) would be also forward 
invariant: \( f(A) = A \) and \( f(B) = B \). Thus no point 
originating in \( A \) could ever intersect \( B \) and vice-versa 
and 
we essentially have two independent dynamical systems. 

Simple examples such as the Dirac-\( \delta_{p} \) measure on a fixed 
point \( p \) are easily shown to be ergodic, but in
general this is a highly non-trivial property to prove. A lot of the 
techniques and methods to be described below are fundamentally 
motivated by the basic question of whether some relevant invariant 
measures are ergodic. It is known that Lebesgue measure is ergodic 
for 
circle rotations \( f(\theta) = \theta + \alpha \) when \( 
\alpha \) is irrational and for covering maps \( f(\theta) = \kappa 
\theta \) when \( \kappa \in \mathbb N \) is \( \geq 2 \) (the proof 
of ergodicity for the latter case will be sketched below). Irrational 
circle rotations are very special because they do not admit any other 
invariant measures besides Lebesgue measure. On the 
other hand covering maps have infinitely many periodic points and 
thus admit infinitely many invariant measures. It is sometimes easier 
to show that certain examples are not ergodic. This is clearly true 
for example for Lebesgue measure and the identity map since any 
subset is fully invariant. A less trivial example is the map 
\( f: [0,1]\to [0,1] \) given by 
 \[ 
 f(x) = \begin{cases}
 2x &\text{ if } 0\leq x \leq 1/4 \\
 -2x+1 &\text{ if } 1/4 \leq x < 1/2 \\
 2x-1/2 &\text{ if } 1/2\leq x \leq 3/4\\
 -2x + 5/2 &\text{ if } 3/4\leq x \leq 1.
 \end{cases}
 \]
Lebesgue measure is invariant , but
 the intervals \( [0, 1/2)  \) and \( [1/2, 1] \) are both 
 backward (and forward) invariant.  This 
 example can easily be generalized by by defining two different 
 Lebesgue measure preserving transformations mapping each of the two 
 subintervals  \( [0, 1/2)  \) and \( [1/2, 1] \)  into themselves. 

The fundamental role played by the notion of ergodicity is given by 
the well known and classical \emph{Birkhoff Ergodic Theorem}. 
We give here only a special case of this result. 
\index{Birkhoff's Ergodic Theorem}
\begin{theorem*}[\cites{Bir31, Bir42}] 
    Let \( f: M \to M \) be a measurable map and 
    let \( \mu \) be an ergodic invariant probability measure for \( 
f \).  
   Then, for any function \( \varphi: M \to \mathbb R \) in \( 
   \mathcal L^{1}(\mu) \), i.e. such that \( \int \varphi d\mu < 
\infty 
   \), and for \( \mu \) almost every \( x \) we have
   \[ 
  \frac{1}{n}\sum_{i=1}^{n}\varphi (f^{i}(x)) \to \int \varphi d\mu
   \]
In particular, for any measurable set \( A\subset M \), letting \( 
\varphi=\mathbb 1_{A} \) be the characteristic function of \( A \), 
we 
have for \( \mu \) almost 
   every \( x\in M \),
 \begin{equation}\label{birkhoff}
\frac{\#\{1\leq j \leq n : f^{j}(x) \in A\}}{n} \to \mu(A).
 \end{equation}
 \end{theorem*}
Here \( \#\{1\leq j \leq n : f^{j}(x) \in A\} \) denotes the
cardinality of the set of indices \( j \) for which \( f^{j}(x) \in A 
\).  
Thus the \emph{average proportion of time} 
which the orbit of a typical point spends in \( A \) 
converges \emph{precisely} to the \( \mu \)-measure of \( A \). 
Notice that 
the convergence of this proportion as \( n\to\infty \) is in itself 
an extremely remarkable and non-intuitive result. The fact that the 
limit is given \emph{a priori} by \( \mu(A) \) means in particular 
that this limit is \emph{independent of the specific initial 
condition} \( x \).  Thus \( \mu \)-almost every initial condition 
has 
the same statistical distribution in space and this distribution 
depends only on \( \mu \) and not even on the map \( f \), except 
implicitly for the fact that \( \mu \) is ergodic and invariant for 
\( f \).

\subsection{Absolute continuity} 

Some care needs to be taken when applying Birkhoff's ergodic theorem 
to 
maps which admit several ergodic invariant measures. Consider for 
example the circle map \( f(\theta) = 10 \theta \). This maps 
preserves Lebesgue measure and also has several fixed points, e.g. 
\( p=0.2222\ldots \), on which we can consider the Dirac-\( 
\delta_{p} 
\) measure.  Both these measures are ergodic. Thus an application of 
Birkhoff's theorem says that ``almost every'' point spends an average 
proportion of time converging to \( m(A) \) in the set \( A \) but 
also that ``almost every'' point spends an average proportion of time 
converging to \( \delta_{p}(A) \) in the set \( A \). If \( m(A) 
\neq \delta_{p}(A) \) this may appear to generate a contradiction. 

The crucial observation here is that the notion of \emph{almost 
every} point is always understood \emph{with respect to a particular 
measure}. Thus Birkhoff's ergodic theorem asserts that for a given 
measure \( \mu \) there exists a set \( \tilde M \subset M \) with \( 
\mu(\tilde M) = 1 \) such that the convergence property holds for 
every \( 
x\in\tilde M \) and in general it may not be possible to identify \( 
\tilde M \) explicitly. Conversely, if \( X\in M \) satisfies \( \mu 
(X) = 0 \) then no conclusion can be drawn about whether 
\eqref{birkhoff} holds for any point of \( X \). 
Returning to the example given above we have the following situation: 
the convergence \eqref{birkhoff} of the time averages to \( 
\delta_{p}(A) 
\) can be guaranteed only for points belonging to a minimal set of 
full measure. But in this case this set reduces to the single point 
\( 
p \) for which \eqref{birkhoff} clearly holds. On the other hand the 
single point \( p \) clearly has zero Lebesgue measure and thus the 
convergence
\eqref{birkhoff} to \( m(A) \) is not guaranteed by Birkhoff's 
Theorem. Thus there is no contradiction. 

An important point therefore is that the information provided by 
Birkhoff's ergodic theorem depends on the measure \( \mu \) 
under consideration. Based on the premise that Lebesgue is the given 
``physical'' measure and that we consider a satisfactory description 
of the dynamics one which accounts for a sufficiently large set of 
points from the point of view of Lebesgue measure, it is clear that 
if \( \mu \) is a Dirac-\( \delta \) measure on a 
fixed point it gives essentially no useful information. On the other 
hand, if \( \mu \) is Lebesgue measure itself then we do get a 
convergence result that holds for Lebesgue almost every starting 
condition.  The invariance of Lebesgue measure is a very special 
property but much more generally we can ask about the existence of 
ergodic invariant measures \( \mu \) which are 
\emph{absolutely continuous} with respect to \( m \). 
\index{Absolutely continuous measure}
\begin{definition}
\( \mu \) is absolutely continuous with respect to \( m \) if 
\[
m(A) = 0\quad\text{ implies} \quad  \mu(A)=0 
 \]
for every measurable set \( A\subset M \). 
\end{definition}
In this case, Birkhoff's theorem implies that \eqref{birkhoff} holds 
for 
all points belonging to a set \( \tilde M\subset M \) with \( 
\mu(\tilde M) = 1 \) and the absolute continuity of \( \mu \) with 
respect to \( m \) therefore implies that \( m(\tilde M) >0 \). 
\index{acip}
Thus the existence of an ergodic \emph{absolutely continuous 
invariant 
probability (acip)} \( \mu \) implies some control over the 
asymptotic 
distribution of at least a set of positive Lebesgue measure. It also 
implies that such points tend to have a dynamics which is non-trivial 
in the sense that it is distributed 
over some relatively large subset of the space as opposed to 
converging for example to some attracting fixed point. Thus it 
indicates that there is a minimum amount of inherent 
\emph{complexity} 
as well as \emph{structure}.  
This 
motivates the basic question: 

\medskip
\begin{center}
\shadowbox{ 
1) Under what conditions does  \( f \) admit an ergodic acip? 
}
\end{center}

This question is already addressed explicitly by Hopf \cite{Hop32}
for invertible transformations. Interestingly he formulates some 
conditions in terms of the existence of what are essentially some 
induced transformations, similar in some respects to the Markov 
structures to be defined below. \footnote{Hopf's result is the 
following: suppose that for every measurable partition \( \mathcal P 
\) 
of the manifold \( M \) and every stopping time function \( p \) 
such that the images 
\( f^{p(\omega)}(\omega) \) for \( \omega\in\mathcal P \) 
are all disjoint, the union of all 
images has full measure. Then \( f \) admits an absolutely continuous 
invariant probability measure. }
In these notes we shall discuss what is effectively a generalization 
of this basic approach.

 \subsection{Mixing}  
\index{mixing}
 Birkhoff's ergodic theorem is very powerful but it is easy to see 
 that the asymptotic space distribution given by \eqref{birkhoff} 
does 
 not necessarily tell the whole story about the dynamics of a given 
 map \( f \). Indeed these conclusions depend not on \( f \) but 
 simply on the fact that Lebesgue measure is invariant and ergodic. 
 Thus from this point of view  the dynamics of an irrational 
 circle rotation \( f(\theta) = \theta + \alpha \) and of the map \( 
 f(\theta) = 2\theta \) are indistinguishable. However it is clear 
 that they give rise to very different kinds of dynamics. In one case 
 for example, nearby points remain nearby for all time, whereas in 
  the other they tend to move apart at an exponential speed. This 
  creates a kind of \emph{unpredictability} in one case which is not 
  present in the other.
  
  \begin{definition}
We say that   
    an invariant probability measure \( \mu \) is \emph{mixing} if 
    \[ 
    |\mu(A\cap f^{-n}(B))-\mu(A)\mu(B)| \to 0 
    \]
    as \( n\to \infty \), for all measurable sets \( A, B \subseteq M 
\). 
\end{definition}
    Notice that mixing implies ergodicity and is therefore a stronger 
property. Thus a natural follow up to 
question \( 1 \) is the following. Suppose that \( f \) admits an 
ergodic \emph{acip} \( \mu \). 

\medskip
\begin{center}
\shadowbox{2) Under what conditions is \( \mu \) mixing ?}
\end{center}

Early work in ergodic 
theory in the 1940's 
considered the question of the genericity of the mixing property is 
in various spaces of systems
\cites{Hal44, Roh48, Ros56, Hol57, 
KacKes58} but, as with ergodicity, in specific classes of systems 
 it is generally easier to show that a system is 
not mixing rather than that it is mixing. For example it is 
immediate that irrational circle rotations are not mixing. On the 
other hand it is non-trivial that 
maps of the form \( f(\theta) = \kappa \theta \) for 
integers \(  \kappa \geq 2 \)  are mixing. 

To develop an intuition for the concept of mixing, notice that 
mixing  
is  equivalent to the condition 
  \[ 
  \left|\frac{\mu(A\cap f^{-n}(B))}{\mu(B)} -\mu(A)\right| \to 0
 \] 
as  \( n\to \infty \), for all measurable sets 
\( A, B \subseteq M \), with \( \mu(B) 
    \neq 0 \). 
 In this form there are two natural interpretations of mixing, one 
 geometrical and one probabilistic. From a geometrical point of view, 
 recall that  \( \mu(f^{-n}(B)) = \mu(B) \) by the invariance of the 
measure.
Then one can think of
 \( f^{-n}(B) \) as a ``redistribution of mass'' and 
 the mixing condition says that for large \( n \) the proportion
 of  \( f^{-n}(B) \) which intersects \( A \) is just proportional to 
the
 measure of \( A \).  In other words \( f^{-n}(B) \) is spreading
 itself uniformly with respect to the measure \( \mu \).  
A more probabilistic point of view is to think of \( \mu(A\cap
f^{-n}(B))/\mu(B) \) as the conditional probability of having \(
x \in A \) given that \( f^{n}(x)\in B \), i.e. the
probability that the occurrence of the event \( B \) today is a
consequence of the occurrence of the event \( A \) \( n \)
steps in the past.  The mixing condition then says that this 
probability converges  to the probability of \( A \), i.e., 
asymptotically,
there is no causal relation between the two events. 
This is why we say that a mixing system exhibits 
\emph{stochastic-like} 
or \emph{random-like} behaviour.

\subsection{Decay of correlations}
It turns out that mixing is indeed a quite generic property at least 
under certain assumptions which will generally hold in the examples 
we shall be interested in. Thus apparently very different systems 
admit mixing \emph{acip}'s and become, in some sense, statistically 
indistinguishable at this level of description. 
Thus it is natural to want to dig deeper in an attempt relate finer 
statistical properties with specific geometric characteristics of 
systems under considerations. 
One 
way to do this is to try to distinguish systems which mix at 
different \emph{speeds}. To formalize this idea we need to generalize 
the definition of mixing. Notice first of all that the original 
definition can be written in integral form as
 \[ 
 \left|\int \mathbb 1_{A\cap f^{-n}(B)} d\mu - \int \mathbb 1_{A} 
d\mu \int
 \mathbb 1_{B} d\mu\right| \to 0
 \]
 where \( \mathbb 1_{X} \) denotes the characteristic function of the 
set \( 
 X \). 
This can be written in the equivalent form
 \[ 
  \left|\int \mathbb 1_{A} (\mathbb 1_{B}\circ f^{n}) d\mu - \int 
\mathbb 1_{A} d\mu \int
  \mathbb 1_{B} d\mu\right| \to 0
 \]
 and this last formulation now admits a natural generalization by 
replacing
 the characteristic functions with arbitrary measurable functions.
 \index{Decay of correlations}
 \begin{definition}
     For real valued measurable functions \( \varphi, \psi : M \to 
\mathbb
     R \) we define the \emph{correlation function}
     \footnote{The derivation of the correlation 
     function from the definition of mixing as given here does not 
perhaps 
     correspond to the historical development. I believe that the 
notion of decay of 
     correlation arose in the context of statistical mechanics and 
was 
     not directly linked to abstract dynamical systems framework 
until 
     the work of Bowen, Lebowitz, Ruelle and Sinai in the 1960's and 
1970's. 
     \cites{Sin68, Bow70, Sin72, PenLeb74, Bow75}  }
\[ 
     \mathcal C_{n}(\varphi, \psi) = \left|\int \psi (\varphi\circ 
f^{n})
     d\mu - \int \psi d\mu \int \varphi d\mu\right|
\]
\end{definition} 
In this context, the functions \( \varphi \) and \( \psi \) are often
 called \emph{observables}.  If \( \mu \) is mixing, the correlation 
 function decays to zero whenever the observables \( \phi, \psi \) 
 are characteristic functions. It is possible to show that indeed it 
 decays also for many other classes of functions. We then have the 
 following very natural question. 
Suppose that the measure \( \mu \) is mixing, fix two observables \( 
\varphi, \psi,\)  and let \( \mathcal C_{n}=\mathcal C_{n}(\varphi, 
\psi)
\). 
 
 \medskip
 \begin{center}
     \shadowbox{3) Does \( \mathcal C_{n} \) decay at a specific 
\emph{rate} 
     depending only on \( f \) ?}
  \end{center}
  
The idea behind this question is that a system may have an
intrinsic \emph{rate of mixing} which reflects some characteristic 
geometrical structures. It turns out that an intrinsic rate does 
sometimes exist and is in some cases possible to determine, but only 
by 
restricting to a suitable class of observables. Indeed, a classical 
result says that even in the ``best'' cases it is possible to choose 
subsets \( A, B \) such that the correlation function \( \mathcal 
C_{n}(\mathbb 1_{A}, \mathbb 1_{B}) \) of the corresponding 
characteristic functions decays at an arbitrarily slow rate. Instead 
positive results exist in many cases by restricting to, for example, 
the space of observables of bounded variation, or H\"older 
continuous, 
or even continuous with non-H\"older modulus of continuity. 
Once the space \( \mathcal H \) of observables has been fixed, the 
goal is to show 
that there exists  a sequence \( \gamma_{n}\to 0 \) (e.g. \(
 \gamma_{n}=e^{-\alpha n} \) or \( \gamma_{n}=n^{-\alpha} \) for some 
 \( \alpha> 0 \))
 depending only on \( f \) and \( \mathcal H \), 
 such that for any two \( \varphi, \psi \in\mathcal H\)
 there exists a constant \( C=C(\varphi, \psi) \) (generally 
 depending on the observables \( \varphi, \psi \)) such that
 \[ 
 \mathcal C_{n} \leq C \gamma_{n}
 \]
 for all \( n\geq 1 \). 
Ideally we would like to show that \( \mathcal C_{n} \) actually 
decays like \( \gamma_{n} \), i.e. to have both lower and uppoer 
bounds, but this is known only in some very 
particular cases. Most known results at present are upper bounds and 
thus when we say that the correlation functions decays at a certain 
rate we will usually mean that it decays \emph{at least} at that 
rate. 
Also, most known results deal with H\"older continuous observables 
and 
thus, to simplify the presentation, we shall assume assume that we 
are 
dealing with this class unless we mention otherwise. 

We shall discuss below several examples of 
systems whose correlations decay at different rates, for example 
\emph{exponential, polynomial} or even \emph{logarithmic}, and a 
basic 
theme of these notes will be gain some understanding about \emph{how} 
and \emph{why} such differences occur and what this tells us about 
the system.

\section{Markov Structures}
\label{markov}

\index{Markov map}
\begin{definition}
    \( f: M \to M \) is (or admits) a \emph{Markov map} if 
there exists a finite or countable  partition \( \mathcal P \) (mod 
0) of \( M \) into 
open sets with smooth boundaries such that \( f(\omega) = M  \)
for every partition element \( 
\omega \in \mathcal P \) and  \( f|_{\omega} \) is a 
continuous non-singular bijection.  
\end{definition}
We recall that a partition mod 0 of \( M \) means that Lebesgue 
almost 
every point belongs to the interior of some partition elements. 
Also, \( f|_{\omega} \) is non-singular if \( |A| > 0 \) implies \( 
|f(A)| > 0 \) for every (measurable) \( A\subset \omega \). 
These two conditions together immediately imply that the full forward 
orbit of almost every point always lies in the interior of some 
partition element. 
The condition \( f(\omega) = M \) 
is a particularly strong version of what is generally referred to as 
the Markov property where it is only required that the image of each 
\( 
\omega \) be a continuous non-singular bijection onto some union of 
partition elements and not necessarily all of \( M \). The stronger 
requirement we use here is sometimes called a Bernoulli property. 
A significant generalization of this definition allows the partition 
element to be just measurable sets and not necessarily open; the 
general results to be given below apply in this case also. However we 
shall not need this for any of the applications which we shall 
discuss.

A natural but extremely far-reaching generalization of the notion of 
a Markov map is the following. 
\index{Induced Markov map}
\index{Markov induced map}
\begin{definition}
$f: M \to M $ admits an \emph{induced Markov map} if 
there exists an open set $ \Delta \subset M$, 
a partition \( \mathcal P \) (mod 0) of 
$\Delta$ and a return time
function $R:\Delta \to \mathbb{N}$, piecewise 
constant on each element of  $\mathcal P$,
such that  the induced map 
\(
F:\Delta \to \Delta 
\)
 defined by
 \(
 F(x) = f^{R(x)}
\)
is a Markov map. 
\end{definition}
Again, the condition that \( \Delta \) is open is not strictly 
necessary. 
For the rest of this section  we shall 
suppose that \( F: \Delta \to \Delta \) is an induced 
Markov map associated to some map \( f: M \to M \). We call \( 
\mathcal P \) the Markov partition associated to the Markov map \( F 
\). 
Clearly if \( f \) is a Markov map to begin with,
it trivially admits an induced Markov map 
with \( \Delta= M \) and \( R\equiv 1 \).  
Since \( \mathcal P \) is assumed to be countable, we can define an 
indexing set 
 \( \mathcal \omega = \{0, 1, 2, \ldots\} \) of the 
Markov partition \( \mathcal P \).  Then, for any finite sequence   
\( a_{0}a_{1}a_{2}n\ldots, a_{n} \) with \( a_{i}\in \mathcal I \), 
we can define the 
\emph{cylinder set of order} \( n \) by 
\[
\omega^{(n)}_{a_0a_1\ldots a_{n}}  \left\{ x : F^i(x)\in \omega_{a_i} 
\text{ 
for } 0 \leq i 
\leq  n\right\}.
\]
Inductively, given $\omega_{a_0a_1\ldots a_{n-1}}$, then 
$\omega_{a_0a_1\ldots
a_n}$ is the part of $\omega_{a_0a_1\ldots a_{n-1}}$ mapped to 
$\omega_{a_n}$ 
by $F^n$.
The cylinder sets define refinements of the partition $\mathcal P$. 
We let \( \omega^{(0)} \) denote generic elements of \( \mathcal 
P^{(0)}= 
\mathcal P \) and \( \omega^{(n)}  \) denote generic elements  of \( 
\mathcal P^{(n)} \). Notice that by the non-singularity of the map \( 
F \) on each partition element and the fact that \( \mathcal P \) is 
a 
partition mod 0, it follows that each \( \mathcal P^{(n)} \) is also 
a 
partition mod 0 and that Lebesgue almost every point in \( 
\Delta \) falls in the interior of some partition element of \( 
\mathcal P \) for all future iterates. In particular almost every \( 
x\in \Delta \) has an associated infinite symbolic sequence \( 
\underline a(x) \) 
determined by the future iterates of \( x \) in relation to the 
partition \( \mathcal P \). 
To  get the much more sophisticated results on the 
statistical properties of \( f \) we need first of all the following 
two additional conditions. 

\begin{definition}
 \( F: \Delta \to \Delta \) has \emph{integrable} or 
\emph{summable} return times if 
\[ 
\int_{\Delta} R(x) dx = \sum_{\omega\in P} |\omega| R(\omega) < 
\infty.
\]
\end{definition}
\index{Bounded distortion}
\begin{definition} \label{def:dist}
\( F: \Delta\to \Delta \)  has the \emph{geometric
   self-similarity} property (or, more prosaically, the (volume) 
\emph{bounded 
   distortion} property) 
   if there exists a constant \( \mathcal D>0 \) such that for 
all \( n\geq 1 \) and any measurable subset \( 
\tilde\omega^{(n)}\subset 
\omega^{(n)}\in\mathcal P^{(n)} \) we have 
     \begin{equation}\label{eq:dist}
\frac{1}{\mathcal D} \frac{|\tilde\omega^{(n)}|}{|\omega^{(n)}|}  
\leq   
\frac{|f^n(\tilde\omega^{(n)})|}{|f^n(\omega^{(n)})|}
\leq  \mathcal D \frac{|\tilde\omega^{(n)}|}{|\omega^{(n)}|}.
    \end{equation}
  \end{definition}

 This means that the relative measure of subsets of a
  cylinder set of any level \( n \) are preserved up to some factor 
\( \mathcal D \) under iteration by \( f^{n} \).  A crucial 
observation here is that the constant \( \mathcal D \) is independent 
of \( n \). Thus in some sense the geometrical structure of any 
subset of \( \Delta \) re-occurs at every scale inside each partition 
element of \( \mathcal P^{(n)} \) up to some bounded distortion 
factor. 
This is in principle a very strong condition but we shall se below 
that it is possible to verify it in many situations. 
We shall discuss in the next section 
some techniques for verifying this 
condition in practice.  First of all we state the first result of 
this 
section. 

\begin{theorem}\label{MarkovMaps}
  Suppose that \( f: M  \to M \) admits an induced 
  Markov map satisfying the  geometric self-similarity
  property and having summable return times. Then it admits an
  ergodic absolutely continuous invariant probability measure \( \mu 
\).
\end{theorem} 

This result goes back to the 1950's 
and is often referred to as the \emph{Folklore Theorem} of 
dynamics. 
We will sketch below the main ideas of the proof. 
First however we state a much more recent 
result which applies in the same setting 
but takes the conclusions much further in the direction of mixing and 
rates of decay of correlations. 
First of all we shall assume without loss of generality that the 
greatest common divisor of all 
values taken by the return time function \( R \) is 1. If this were 
not the case all 
return times would be multiples of some integer \( k\geq 2 \) and the 
measure \( \mu \) given by 
the Theorem stated above would clearly not be mixing. If this is the 
case however, we could just consider the map \( \tilde f = f^{k} \) 
and the results to be stated below will apply to \( \tilde f \) 
instead of \( f \). 
We define the 
\emph{tail of the return times} as the measure of the set 
\[ 
R_{n} = \{x\in\Delta : R(x) > n \}
\]
of points whose return times is strictly larger than \( n \). 
The integrability condition implies that 
\( R(x) < \infty \) for almost every point 
and thus 
\[
|R_{n}|\to 0 
\]
as \( n\to\infty \). 
However there is a range of possible \emph{rates} of decay of \( 
|R_{n}| \) all of which are compatible with the integrability 
condition. 
L.-S. Young observed and proved that a bound on the decay of 
correlations for H\"older continuous observables can be obtained from 
bounds on the rate of decay of the tail of the return times. 

\begin{theorem}[\cites{You98,You99}]\label{MarkovTowers}
    Suppose that \( f: M  \to M \) admits an induced 
    Markov map satisfying the geometric self-similarity
    property and having summable return times. Then it admits an
    ergodic (and mixing) 
    absolutely continuous invariant probability measure.
Moreover  the correlation function for H\"older continuous 
observables satisfies the following bounds: 
\begin{description}
    \item[Exponential tail] 
    If \( \exists  \)\( \alpha > 0 \) such that 
    \( 
    |R_{n}|=\mathcal O(e^{\alpha n}) \), then \( \exists  \) \( 
    \tilde\alpha > 0 \) such that \( \mathcal C_{n} = \mathcal 
    O(e^{-\tilde\alpha n}) \).
    \item[Polynomial tail]
If \( \exists \)\( \alpha > 1 \) such that \( 
|R_{n}|=\mathcal O(n^{\alpha}) \), then \( \mathcal C_{n} = \mathcal 
O(n^{-\alpha+1}) \). 
\end{description}
\end{theorem}
 Other papers have also addressed the question of the decay of 
 correlations for similar setups mainly using spectral operator 
 methods 
\cites{You98, Bre99, Mau01a,  BreFerGal99, Mau01b}. 
We remark that the results about the rates of decay of correlations 
generally require an a priori slightly stronger form of bounded 
distortion than that given in \eqref{eq:dist}. The proof 
in \cite{You99} uses a very geometrical/probabilistic 
\emph{coupling argument} which appears to be quite versatile and 
flexible. Variations of the argument have been applied to prove the 
following generalizations which apply in 
the same setting as above (in both cases we state only a 
particular case of the theorems proved in the cited papers). 

The first one extends Young's result to 
arbitrarily slow rates of decay. 
We say that  \( \rho: \mathbb R^{+}\to\mathbb R^{+} 
\) is \emph{slowly varying} (see \cite{Aar97})  if for all \( y> 0 \) 
we have \( 
\lim_{x\to\infty} \rho (xy)/\rho(x) = 0 \).  A simple example of a 
slowly varying function is the function \( \rho(x) = e^{(\log 
x)/(\log\log x)} \). Let \( \hat R_{n}=\sum_{\hat n \geq n} R_{\hat 
n}\). 

\begin{theorem}[Hol04]\label{slowdecay}
    The correlation function for H\"older 
    continuous observables satisfies the following bound. 
\begin{description}
    \item[Slowly varying tail] If \( \hat R_{n}=\mathcal 
    O(\rho (n)) \) where \( \rho \) is a monotonically decreasing to 
    zero, slowly varying, \( C^{\infty} \) function, then 
   \( \mathcal C_{n}=\mathcal  O(\rho(n)) \). 
 \end{description}   
\end{theorem}

The second extends Young's result  to observables with 
very weak, non-H\"older, modulus of continuity. 
We say that  \( \psi : I \to \mathbb R\) has a \emph{logarithmic
modulus of continuity \( \gamma \)} if there exists \( C > 0\) such 
that for 
all \( x, y \in I \) we have  
    \[ 
   |\psi (x) - \psi(y) | \leq C | \log |x-y||^{-\gamma}.
    \]
 For both the exponential and 
 polynomial tail situations we have the following   
\begin{theorem}[Lyn04]\label{nonholder}
There exists  \( \alpha>0 \) such that for all \( \gamma \) 
    sufficiently large and  observables with logarithmic modulus of 
    continuity \( \gamma \), we have \( \mathcal C_{n} =\mathcal 
    O(n^{-\alpha})\).
  \end{theorem}

These general results indicate that the rate of decay of correlations 
is linked to what is in effect the \emph{geometrical structure} of \( 
f 
\) as reflected in the tail of the return times for the induced map 
\( 
F \).  From a technical point of view they shift the problem of the 
statistical properties of \( f \) to the problem of the 
geometrical structure of \( f \) and thus to the (still highly 
non-trivial) problem of showing that \( f \) admits an induced Markov 
map and of estimating the tail of the return times of this map. 
The construction of an induced map 
in certain examples is relatively straightforward and essentially 
canonical but the most interesting constructions require statistical 
arguments to even show that such a map exists and to estimate 
the tail of the return times. In these cases the construction is 
not canonical and it is usually not completely clear to what 
extent the estimates might depend on the construction. 
  
 We now give a sketch of the  proof of Theorem \ref{MarkovMaps}. 
 The proofs of Theorems \ref{MarkovTowers}, \ref{slowdecay} and 
 \ref{nonholder} are in a similar spirit and we refer the interested 
 reader to the original papers. We assume throughout the next few 
 sections that \( F: \Delta \to \Delta \) is the Markov 
induced map associated to \(  f: I \to I\) and \( \mathcal P^{(n)} \) 
are the family of cylinder sets generated by the Markov partition \( 
\mathcal P=\mathcal P^{(0)} \) of \( \Delta \). 
We first define a measure \( 
\nu \) on \( \Delta \) and show in  that it is \( F 
\)-invariant, ergodic, and absolutely continuous with respect to 
Lebesgue. Then we define the measure \( \mu \) on \( I \) in terms 
of \( \nu \) and show that it is \( f \)-invariant, ergodic, and 
absolutely continuous.

 \subsection{The invariant measure for \protect\( F \protect\)}
We start with a
preliminary result which is a consequence of the bounded distortion 
property. 
\subsubsection{The measure of cylinder sets}
A straightforward but remarkable consequence of the bounded 
distortion 
property is that the measure of cylinder sets tends to zero 
uniformly. 

\begin{lemma}\label{cylindersets}
\[
\max \{|\omega^{(n)}|; \omega^{(n)}\in\mathcal P^{(n)}\} \to 0 
\quad \text{ as } \quad n\to 0
\]
  \end{lemma} 

  Notice that in the one-dimensional case, the measure of an interval 
coincides with its diameter and so this implies in particular that 
the diameter of cylinder sets tends to zero, implying 
the essential uniqueness of the symbolic representation of 
itineraries. 
 \begin{proof}
     It is sufficient to show that there exists a constant \( \tau 
     \in (0,1) \) such that for every \( n\geq 0  \) and every \( 
     \omega^{(n)}\subset \omega^{(n-1)} \) we have 
     \begin{equation}\label{cylconv}
    |\omega^{(n)}|/|\omega^{(n-1)}| \leq \tau.
     \end{equation}
 Applying this inequality recursively then implies 
     \(
    |\omega^{(n)}| \leq \tau |\omega^{(n-1)}| \leq 
\tau^{2}|\omega^{(n-2)}| \leq 
    \dots \leq \tau^{n} | \omega^{0}| \leq \tau^{n}|\Delta|.
     \)
 To verify 
     \eqref{cylconv}  we shall show that 
     \begin{equation}\label{cylconv2}
     1-\frac{|\omega^{(n)}|}{|\omega^{(n-1)}|} = 
\frac{|\omega^{(n-1)}|-|\omega^{(n)}|}{\omega^{(n)}|}=
\frac{|\omega^{(n-1)}\setminus \omega^{(n)}|}{|\omega^{(n)}|} \geq 
1-\tau.
     \end{equation}
     To prove \eqref{cylconv2} let first of all 
 \(
     \delta = \max_{\omega\in\mathcal 
       P}|\omega| < |\Delta|. 
 \)
Then,  from the definition of cylinder sets we have that 
  \( F^{n}(\omega^{(n-1)}) = \Delta \) and that \( 
F^{n}(\omega^{(n)}) 
  \in \mathcal P = \mathcal P^{(0)} \), and therefore 
  \(  |F^{n}(\omega^{(n)})| \leq \delta \) or, equivalently, 
\(
|F^{n}(\omega^{(n-1)}\setminus \omega^{(n)}|)| \geq |\Delta| - 
  \delta > 0. 
\) 
Thus, using the bounded distortion 
 property we have 
\[
     \frac{|\omega^{(n-1)}\setminus \omega^{(n)}|}{|\omega^{(n)}|}  
     \geq \frac{1}{\mathcal D} 
     \frac{|F^{n}(\omega^{(n-1)}\setminus \omega^{(n)}|)|}
     {|F^{n}(\omega^{(n)})|} 
    \geq \frac{|\Delta| - 
  \delta}{|\Delta| \mathcal D}
\]
and \eqref{cylconv2} follows choosing  \( \tau = 1- ({(|\Delta| - 
  \delta)}/{|\Delta| \mathcal D}) \).
 \end{proof}

The next property actually follows only from the conclusions of Lemma 
\ref{cylindersets} rather than from the bounded distortion property 
itself. 
It essentially says that it is possible to  ``zoom in'' to any 
given set of positive measure.

  \begin{lemma}\label{ergolemma}
For any \( \varepsilon > 0 \) and any Borel set \( A 
      \) with \( |A| > 0 \) there exists \( n\geq 1 \) and \( 
      \omega^{(n)}\in\mathcal P^{(n)} \) such that 
      \[ 
     |A\cap \omega^{(n)}| \geq (1-\varepsilon) |\omega^{(n)}|.
      \]
      \end{lemma}
\begin{proof}
Fix some \( \varepsilon > 0 \). 
Suppose first of all that \( A \) is compact. Then, using the 
properties of Lebesgue measure it is possible to show that 
for any \( \eta>0 \) there exists 
an integer \( n\geq 1 \) and a collection 
\( \mathcal I_{\eta}=\{\omega^{n}\}\subset \mathcal P^{(n)}
     \)
     such that 
     \(
     A \subset \cup_{\mathcal \omega_{\eta} }\omega^{(n)} 
     \)
     and 
     \(
     |\mathcal \omega_{\eta}| \leq |A| + \eta. 
     \) 
 Now suppose by contradiction that 
  \( |\omega^{(n)}\cap A| \leq (1-\varepsilon) |\omega^{(n)}| \) for 
every \( 
     \omega^{(n)}\in \mathcal \omega_{\eta} \) for any given \( \eta 
> 0 \). 
   Using that fact that the \( \omega^{(n)}\in 
    \mathcal \omega_{\eta} \) are disjoint and thus  \( 
    \sum |\omega^{(n)}| = |\mathcal \omega_{\eta}| \), this implies 
that 
\[
	 |A|= \sum_{\omega^{(n)}\in\mathcal \omega_{\eta}}|\omega^{(n)}\cap 
A| 
	 \leq (1-\varepsilon) \sum_{\omega^{(n)}\in\mathcal \omega_{\eta}} 
	 |\omega^{(n)}| \leq (1-\varepsilon) (|A| + \eta).
\]
Since \( \eta \) can be chosen arbitrarily small after fixing \( 
\varepsilon \) this gives a contradiction. 
If A is not compact we can approximate if from below in measure by 
compact sets and repeat essentially the same argument. 
\end{proof}

\subsubsection{Absolute continuity}
The following estimate also follows immediately from the bounded 
distortion property. It says that the absolute continuity property of 
\( 
F \) on partition elements is preserved up to arbitrary scale with 
uniform bounds. 

\begin{lemma}\label{abscon}
Let \( A\subset \Delta\) and \( n\geq 1 \). 
Then 
\[ 
|F^{-n}(A)|  
  \leq  \mathcal D |A|. 
\]
\end{lemma}
\begin{proof}
The 
Markov property implies that 
\( F^{-n}(A) \)  is a union of disjoint sets each  
contained in the interior of some element \( \omega^{(n)}\in \mathcal 
P^{(n)} \). Moreover each \( \omega^{(n)} \) is mapped by \( F^{n} \) 
to \( 
 \Delta \) with uniformly bounded distortion, thus we have 
\(
{|F^{-n}(A)\cap \omega^{(n)}|}/{\omega^{(n)}} \leq \mathcal D 
{|A|}/{|\Delta|}
\)
or, equivalently, 
\(
|F^{-n}(A)\cap \omega^{(n)}| \leq {\mathcal D |A| 
|\omega^{(n)}|}/{|\Delta| }. 
\)
Therefore 
\[
|F^{-n}(A)|   = \sum_{\omega^{(n)}\in\mathcal P^{(n)}} 
    |F^{-n}(A)\cap \omega^{(n)}|  
     \leq \frac{\mathcal D |A|}{|\Delta|} 
\sum_{\omega^{(n)}\in\mathcal P^{(n)}} |\omega^{(n)}| 
  = \mathcal D |A|
\] 

\end{proof}

\subsubsection{The pull-back of a measure}
For any \( n\geq 1 \)  and Borel \( A \subseteq \Delta \), let 
\[ 
\nu_n(A) = \frac{1}{n}\sum_{i=0}^{n} |F^{-i}(A)|. 
\] 
It is easy to see that \( \nu_{n} \) is a probability measure on \( 
\Delta \) and absolutely continuous with respect to Lebesgue. 
Moreover, lemma \ref{abscon} implies that that the absolute 
continuity 
property is uniform in \( n \) and \( A \) in the sense that 
\(
\nu_{n}(A) \leq \mathcal D \ |A|
\)
for any \( A \) and for any \( n\geq 1 \). 
By some standard results of functional analysis, this implies the 
following 
\begin{lemma}
There exists a probability measure \( \nu \) 
and a subsequence \( \{\nu_{n_{k}}\} \) such that, 
for every measurable set \( A \),
\begin{equation}\label{measureconv}
\nu_{n_{k}}(A) \to \nu(A) \leq \mathcal D |A|.
\end{equation}
\end{lemma}
In particular \( A \) is absolutely continuous with respect to 
Lebesgue.

\subsubsection{Invariance}
To show that \( \nu \) is \( F \)-invariant, let \( A\subset \Delta 
\) be a 
measurable set. Then, by \eqref{measureconv} 
we have 
\begin{eqnarray*}
\nu\left(F^{-1}(A)\right)&=&\lim_{k\rightarrow\infty} 
\frac{1}{n_k}\sum_{i=0}^{n_k-1}
|F^{-(i+1)}(A)|\\
&=&\lim_{k\rightarrow\infty}\left[\frac{1}{n_k}\sum_{i=0}^{n_k-1}
|F^{-1}(A)|
-\frac{|A|}{n_k} + \frac{|F^{-n_k}(A)|}{n_k}
\right]
\\
\end{eqnarray*}
Since \(|A|\) and \( |f^{-n}(A)|\) are both uniformly 
bounded by 1, we have 
\({|A|}/{n_k} \to 0   \) and \( |F^{-n_k}(A)|/{n_k}\to 0 \) as \( 
k\to 
0 \). Therefore 
\begin{align*} 
\nu(F^{-1}(A))&= 
\lim_{k\rightarrow\infty}\left[\frac{1}{n_k}\sum_{i=0}^{n_k-1}
|F^{-1}(A)|
-\frac{|A|}{n_k} + \frac{|F^{-n_k}(A)|}{n_k}
\right] 
\\ &= 
\lim_{k\rightarrow\infty}\frac{1}{n_k}\sum_{i=0}^{n_k-1}
|F^{-1}(A)| = \nu(A).
\end{align*}
Therefore \( \nu \) is \( F \)-invariant. 

\subsubsection{Ergodicity and uniqueness}

Let \( A\subset I \) be a measurable set with \( F^{-1}(A) = A \) and 
\( 
\mu(A) > 0 \). We shall show that \( \nu(A) = |A|= 1 \). This implies 
both 
ergodicity and uniqueness of \( \nu \). Indeed, if \( \tilde \nu \) 
were another such 
measure invariant absolutely continuous measure, 
there would be have to be a set \( B \) with \( F^{-1}(B) = B 
\) and \( \tilde\nu(B) = 1 \). But  in this case we would have also 
\( 
|B| = 1 \) and thus \( A = B \mod 0 \). This is impossible since two 
absolutely continuous invariant measures must have disjoint support. 

To prove that \( |A| = 1 \), let \( 
A^{c}=\Delta\setminus A \) denote the complement of \( A \). Notice 
that 
\( x\in A^{c} \) if and 
only if \( F(x) \in A^{c} \) and therefore \( F(A^{c})=A^{c} \).  
By Lemma
\ref{ergolemma}, for any \( \varepsilon > 0 \) 
there exists some \( n\geq 1 \)  and 
\( \omega^{(n)}\in\mathcal P^{(n)} 
\)  such that 
\( 
|A\cap \omega^{(n)}| \geq (1-\varepsilon) |\omega^{(n)}| 
\)
and therefore 
\[ 
|A^{c}\cap \omega^{(n)}| \leq \varepsilon |\omega^{(n)}|. 
\]
Using that fact that \( F^{n}(\omega^{(n)}) = I \) and the invariance 
of \( 
A^{c} \) have \( F^{n}(\omega^{(n)}\cap A^{c}) = A^{c}  \). The 
bounded distortion property then gives
\[ 
|A^{c}| = \frac{|F^{n}(\omega^{(n)}\cap 
A^{c})|}{|F^{n}(\omega^{(n)})|} \leq 
\mathcal D \frac{|\omega^{(n)}\cap A^{c}|}{|\omega^{(n)}|} \leq 
\mathcal D 
\varepsilon.
\]
Since \( \varepsilon \) is arbitrary this implies \( |A^{c}| = 0 \) 
and thus \( |A| = 1 \).  

\subsection{The invariant measure for \protect\( f \protect\)}
We now show how to define a probability measure \( \mu \) which is 
invariant for the original map \( f \) and satisfies all the required 
properties.

\subsubsection{The probability measure \protect\( \mu \protect\)}
We let \( \nu_{\omega} \)
denote the restriction of \( \nu \) to the partition element \( 
\omega\
\in\mathcal P\), i.e. for any measurable set \( A\subset \Delta \) we 
have \( \nu_{\omega}(A) = \nu(A\cap \omega) \). Then \( \nu(A) = 
\sum_{\omega\in\mathcal P}\nu_{\omega}(A) \). 
Then, for any measurable set \( A\subseteq M
\) (we no longer restrict our attention to \( \Delta \)) 
we define 
\[ 
\hat\mu(A) = \sum_{\omega\in\mathcal P} 
\sum_{j=0}^{R(\omega)-1} \nu_{\omega}(f^{-j}(A)).
\]
Notice that this is a sum of non-negative terms and is uniformly 
bounded since 
\begin{align*}
\hat\mu(A)\leq \hat\mu(M) &= \sum_{\omega\in\mathcal P}
\sum_{j=0}^{R(\omega)-1} \nu_{\omega}(f^{-j}(M))
\\ &= \sum_{\omega\in\mathcal P}\sum_{j=0}^{R(\omega)-1} 
\nu_{\omega}(M) = 
\sum_{\omega\in\mathcal P} R(\omega) \nu (\omega) < \infty 
\end{align*}
by the assumption on the summability of the return times. 
Thus it defines a finite measure on \( M \) and from this we define a 
probability measure by normalizing to get
\[ 
\mu (A) = \hat \mu (A) / \hat \mu (M).
\]

\subsubsection{Absolute continuity}
The absolute continuity of \( \mu \) is an almost immediate 
consequence of the definition and the absolute continuity of \( \mu 
\). 
Indeed,   \( |A|=0 \) 
implies \( \nu(A) = 0 \) which implies \( \nu_{\omega}(A) = 0 \) for 
all 
\( I\omega\in \mathcal P \), which therefore implies that we have 
\[ 
\sum_{j=0}^{R(\omega)-1} \nu_{\omega}(f^{-j}(A)) = 0
\]
and therefore \( \mu(A) = 0 \). 

\subsubsection{Invariance}
Recall first of all that  by definition \( f^{R(\omega)}(\omega) = 
\Delta 
\) for any \( \omega\in\mathcal P \). 
Therefore, for any \( A \subset M \) we have 
\[ 
f^{-R(\omega)}(A) \cap \omega = F|_{\omega}^{-1}(A) \cap \omega,
\]
where \( F|_{\omega} ^{-1}\) denotes the inverse of the restriction 
\( F|_{\omega} \) 
of \( F \)  to \( \omega \) (notice that \( f^{-R(\omega)}(A) \cap 
\omega  = \emptyset \) if \( A\cap \Delta = \emptyset \)). 
In particular, using the invariance of \( \nu \) under \( F \), 
this gives 
\[
\sum_{\omega\in\mathcal P}\nu (f^{-R(\omega)}(A) \cap \omega) 
= 
\sum_{\omega\in\mathcal P}\nu(F^{-1}|_{\omega}(A) \cap \omega)
=\nu(F^{-1}(A)) = \nu(A).
\]
Using this equality we get, for any measurable set \( A \subseteq I 
\), 
\begin{align*}
    \mu (f^{-1}(A)) &=
    \sum_{\omega\in\mathcal P} \sum_{j=0}^{R(\omega)-1} 
    \nu_{\omega}(f^{-(j+1)}(A)) \\
    &= \sum_{\omega\in\mathcal P} \sum_{j=0}^{R(\omega)-1} 
    \nu (f^{-(j+1)}(A)\cap \omega) 
    \\& 
    =\sum_{\omega\in\mathcal P}
    \nu\left[ (f^{-1}(A)\cap \omega) + \cdots + (f^{-R(\omega)}(A) 
\cap 
    \omega)\right] \\
    & =
    \sum_{\omega\in\mathcal P} \sum_{j=1}^{R(\omega)-1} \nu 
(f^{-j}(A)\cap 
    \omega)  
    + \sum_{\omega\in\mathcal P} (f^{-R(\omega)}(A) \cap 
\omega) 
\\ &= \sum_{\omega\in\mathcal P} \sum_{j=1}^{R(\omega)-1} \nu 
(f^{-j}(A)\cap \omega)  
    + \nu(A) \\
    & = 
    \sum_{\omega\in\mathcal P} \sum_{j=0}^{R(\omega)-1} \nu 
(f^{-j}(A)\cap 
    \omega)  
    \\ & =\mu (A).
\end{align*}

\subsubsection{Ergodicity and uniqueness}

Ergodicity of \( \mu \) follows immediately from the ergodicity of 
\( \nu \) since every fully invariant set for of positive measure 
must intersect the image of some partition element \( \omega\) and 
therefore must have positive (and therefore full) measure for \( \nu 
\) and therefore must have full measure for \( \mu \). 
Notice however that we can only claim a limited form of uniqueness 
for the measure \( \mu \). Indeed, the support of \( \mu \) is given 
by 
\[ 
supp (\mu) = \bigcup_{\omega\in\mathcal P}
\bigcup_{j=0}^{R_{k}-1} f^{j}(I_{k)}
\]
which is the union of all the images of all partition elements. Then 
\( \mu \) is indeed the unique ergodic absolutely continuous 
invariant measure on this set. However in a completely  abstract 
setting there 
is no way of saying that \( supp (\mu) = M\) nor that there may not 
be other relevant measures  in \( M \setminus supp (\nu) \). 

\subsection{Expansion and distortion estimates}

The application of the abstract results discussed above to specific 
examples involves three main steps: 
\begin{itemize}
    \item Combinatorial construction of the induced map;
    \item Verification of the bounded distortion property;
    \item Estimation of the tail of the return times function and 
    verification of the integrability of the return times.
 \end{itemize}   

 We shall 
discuss some of these step in some detail in relation to 
 some of the specific case as we go through them below. Here we just 
make a few 
remarks concerning the bounded distortion property and in particular 
the crucial role played by \emph{regularity}  and 
\emph{derivative} conditions in 
these calculations.

We begin with a quite general observation which relates the 
geometric self-similarity 
condition  to a property involving the derivative  of \( F \). 
    Let \( F: \Delta \to \Delta \) be a Markov map which is 
continuously 
    differentiable on each element of the partition \( \mathcal P \). 
    We let \( \det DF^{n} \) denote the determinant of the derivative 
    of the map  \( F^{n} \). 
   
\begin{definition} 
    We say 
that \( F \) has \emph{uniformly bounded derivative 
    distortion} if there exists a constant \( \mathcal D>0 \) such 
that for 
    for all $n\geq1$ and $\omega \in  P^{(n)}$ we have
    \begin{equation}\label{ddist1}
    Dist(f^n, \omega) :=  \max_{x,y\in I^{(n)}}\log 
    \frac{\det DF^{n}(x)}{\det DF^{n} (y)}\leq \mathcal D 
    \end{equation}
  \end{definition}
 
Notice that this is just the infinitesimal version of the 
self-similarity bounded distortion property and indeed it is possible 
to show that this condition implies the geometric 
self-similarity property. In the one-dimensional setting and assuming 
 \( J \subset \omega\) to be an open set, this implication follows 
 immediately from the Mean Value Theorem. Indeed, 
in one dimension the determinant of the derivative is just the 
derivative itself.  Thus, the Mean Value 
Theorem implies that there exists \( x\in I\omega \) such that \( 
| Df^{n}(x)| = | DF^{n}(\omega)|/|\omega| \) and \( y\in J  \) such 
that  \( | DF(y)| = |  DF^{n}(J)|/|J| \). Therefore 
\begin{equation} \label{derivdist}
\frac{|\omega|}{|J |}
\frac{|F^{n}(J )|}{|F^{n}(\omega)|} = 
\frac{|F^{n}(J )|/|J |}{|F^{n}(\omega)|/|\omega| } 
= \frac{|DF^{n}(y)|}{|DF^{n}(x)|} \leq \mathcal D.
\end{equation}
To verify \eqref{ddist1}
we use the chain rule to write 
\[ 
\log \frac{|\det DF^{n}(x)|}{|\det DF^{n}(y)|} = 
\log \prod_{i=0}^{n-1}\frac{|\det DF(F^{i}(x))|}{|\det DF(F^{i}(x))|}
= \sum_{i=0}^{n-1}\log \frac{|\det DF(F^{i}(x))|}{|\det 
DF(F^{i}(y))|}.
\]
Now adding and subtracting \(  {|\det DF(F^{i}(y))|}/{|\det 
DF(F^{i}(y))|}\) 
and using that fact that \( \log (1+x) < x \) for 
\( x>0 \) gives 
\begin{align*}
\log \frac{|\det DF(F^{i}(x))|}{|\det DF(F^{i}(y))|} &\leq 
\log \left(\frac{|\det DF(F^{i}(x))- \det DF(F^{i}(y))|}{|\det 
DF(F^{i}(y))|}  + 
1\right) 
\\ &
\leq \frac{|\det DF(F^{i}(x))- \det DF(F^{i}(y))|}{|\det 
DF(F^{i}(y))|}.
\end{align*}
Therefore we have 
\begin{equation}\label{eq:distest}  
    \log \frac{|\det DF^{n}(x)|}{|\det DF^{n}(y)|} 
    \leq \sum_{i=0}^{n-1}
\frac{|\det DF(F^{i}(x))-\det DF(F^{i}(y))|}{|\det DF(F^{i}(y))|}
\end{equation}

The inequality \eqref{eq:distest} gives us the basic tool for 
verifying the required distortion properties in particular examples. 

 \section{Uniformly Expanding Maps}\label{expanding}
In this section we discuss 
maps which are \emph{uniformly} expanding. 
\index{Uniformly expanding maps}
 \subsection{The smooth/Markov case}
   We say that \( f \) is uniformly expanding if there exist 
constants \(
     C, \lambda > 0 \) such that for all \( x\in M \), all \( v\in 
T_{x}M
     \), and all \( n\geq 0 \), we have 
     \[ 
     \|Df^{n}_{x}(v)\|\geq Ce^{\lambda n}\|v\|.
     \]
  We remark once again that this is a special case of the nonuniiform 
  expansivity condition. 
\index{Decay of correlations ! Exponential}
\begin{theorem}
 Let \( f: M \to M \) be \( C^{2} \) uniformly expanding. Then there
     exists a unique acip \( \mu
     \)\cites{Ren57, Gel59, Par60,  Rue68, Ave68, 
    KrzSzl69, Wat70, Las73}. 
 The measure \( \mu \) 
 is mixing and the correlation function decays exponentially fast. 
 \cites{Sin72, Per74, Bow75, Rue76}
\end{theorem}
 
The references given here use a variety of arguments some of which
use  the remarkable observation that uniformly expanding maps 
are intrinsically \emph{Markov} in the strong sense given above, with 
\( 
\Delta = M \),  
a finite number of partition elements and return time \( R\equiv 1 
\) (this is particularly easy to see in the case of one-dimensional 
circle maps \( f: \mathbb S^{1} \to \mathbb S^{1} \)). Thus the main 
issue here is the verification of the distortion condition. 

One way to show this is to show that there is a uniform upper bound 
independent 
of \( n \) for the sum in \eqref{eq:distest} above. 
Indeed, notice first of all that the expansivity condition implies in 
particular that   
\(     |\det DF(F^{i}(y))| \geq Ce^{\lambda i} \geq  C > 0 \) for 
every \( y \), 
and the  \( C^{2} \) regularity condition implies that \( \det Df \) 
is Lipschitz: there exists \( L > 0 \) such that 
\(
|\det DF(F^{i}(x))-\det DF(F^{i}(y))| \leq L |F^{i}(x)-F^{i}(y)|.
\) for all \( x,y\in M \). 
Substituting these inequalities into 
\eqref{eq:distest} we get 
\begin{equation}\label{distest2}
    \sum_{i=0}^{n-1}\frac{|\det DF(F^{i}(x))-\det 
    DF(F^{i}(y)|}{|\det DF(F^{i}(y))|} 
    \leq \frac{L}{C}   \sum_{i=0}^{n-1}
    | F^{i}(x))- F^{i}(y)|
\end{equation}
The next step,and final, step 
uses the expansivity condition as well as, implicitly, the Markov 
property in a crucial 
way. 
Indeed, let \( diam M \) denote the diameter of \( M \), 
i.e. the maximum distance between any two points in \( M \). 
The definition of \( \mathcal P^{(n)} \) implies that \( \omega \) is 
mapped diffeomorphically to \( M \) by \( F^{n} \) and thus 
\[ 
|diam M| \geq  |F^{n}(x)-F^{n}(y)| \geq C 
e^{\lambda(n-i)}|F^{i}(x)-F^{i}(y)|
\]
for every \( i=0,\ldots,n-1 \). Therefore
\begin{equation}\label{distest3}
\sum_{i=0}^{n-1} |F^{i}(x)-F^{i}(y)| 
\leq \frac{diam M}{C}\sum_{i=0}^{n-1}  e^{-\lambda (n-i)}  
\leq  \frac{diam M}{C}\sum_{i=0}^{\infty}  e^{-\lambda i}   
\end{equation}
Substituting back into \eqref{distest2} and \eqref{eq:distest} gives 
a bound for the distortion which is independent of \( n \).

The regularity condition on \( \det DF \) can be 
weakened somewhat but not completely. 
There exist 
examples of one-dimensional circle maps 
\( f: \mathbb S^{1}\to \mathbb S^{1} \) which are \( C^{1} \) 
uniformly expanding (and thus Markov as above) but for which the 
uniqueness of the absolutely continuous invariant measure fails 
\cites{Qua96, CamQua01}, essentially due to the failure of the 
bounded 
distortion calculation. 
On the other hand,  the 
distortion calculation above goes through with minor modifications as 
long as \( \det DF \) is just H\"older continuous.  In some 
situations, 
such as the one-dimensional \emph{Gauss map} \( f(x) = x^{-1} \mod 1 
\) which is Markov but for which 
the derivative \( Df \) is not even H\"older continuous, one can 
compensate by taking advantage of the large derivative. Then it is 
possible to show directly that the right hand side of 
\eqref{eq:distest} is 
uniformly bounded, even though  \eqref{distest2} does not hold.

\subsection{The non-Markov case}

The general (non-Markov) piecewise expanding case is significantly 
more complicated and even the existence of an absolutely continuous 
invariant measure is no longer guaranteed \cites{LasYor73, 
GorSch89, Qua99, Tsu00a, Buz01a}. One possible problem
is that  the images of the discontinuity set 
can be very badly distributed and 
cause havoc with any kind of \emph{structure}. In the Markov case 
this 
does not happen because the set of discontinuities gets mapped to 
itself by definition.  Also the possibility 
of components being \emph{translated} in different directions can 
destroy on a global level 
the local expansiveness given by the derivative. 
Moreover, where results exist for rates of decay of correlations, 
they do 
not always apply to the case of H\"older continuous observables, as 
technical reasons sometimes require that different functions spaces 
be 
considered which are more compatible with the discontinuous nature of 
the maps. We shall not explicitly comment on the particular classes 
of observables considered in each case. 

In the one-dimensional case these problems are somewhat more 
controllable and relatively simple conditions guaranteeing the 
existence of an ergodic invariant probability measure can be 
formulated even in the case of a countable number of domains of 
smoothness of the map. These essentially require that the size of the 
image of all domains on which the map is \( C^{2} \) be strictly 
positive  and that certain conditions on 
the second derivative are satisfied
\cites{LasYor73, Adl73, Bow77, Bow79}. 
 In the higher dimensional case, the situation is considerably more 
 complicated and there are a variety of possible conditions which can 
 be assumed on the discontinuities. The conditions of 
 \cite{LasYor73} were generalized to the two-dimensional context in 
 \cite{Kel79} and then to arbitrary dimensions in \cites{GorBoy89, 
 Buz00a, Tsu01b}. There are also several other papers which 
 prove similar results under various conditions, we mention 
 \cites{Alv00, BuzKel01, BuzPacSch01, Buz01c, Sau00, BuzSar03}.
 In \cites{Buz99a, Cow02} it is shown that conditions 
 sufficient for the existence of a measure are \emph{generic} in 
 a certain sense within the class of piecewise expanding maps. 
 
Estimates for the decay of correlations have been proved for 
non-Markov 
piecewise smooth maps, although again the techniques have had to be 
considerably generalized. In terms of setting up the basic arguments 
and techniques, a similar role to that played by 
\cite{LasYor73} for the existence of absolutely continuous invariant 
measures can be attributed to \cites{Kel80, HofKel82, Ryc83} 
for the problem of decay of correlations in the one-dimensional 
context. More recently, alternative approaches have been proposed and 
implemented in \cites{Liv95a, Liv95, You98}. The approach of 
\cite{You98} 
has proved particularly suitable for handling some higher dimensional 
cases such as 
\cite{BuzMau02} in which assumptions on the discontinuity set 
are formulated in terms of \emph{topological pressure}
and \cites{AlvLuzPin, Gou04} in which 
they are formulated  as \emph{geometrical non-degeneracy} 
assumptions and \emph{dynamical } 
assumptions on the \emph{rate of recurrence} 
of typical points to the discontinuities. 
The construction of an induced Markov map is combined in 
\cite{Dia04} with the Theorem \ref{nonholder} to obtain estimates for 
the decay of correlations of non-H\"older observables for Lorenz-like 
expanding maps. 
We remark also that the results of \cites{AlvLuzPin, Gou04} 
apply to more general piecewise nonuniformly expanding maps, see 
section \ref{general}. 
It would be interesting to  
understand the relation between
the assumptions of \cites{AlvLuzPin, Gou04} and 
those of \cite{BuzMau02}. 

 \section{Almost Uniformly Expanding Maps}
  \label{almost expanding}
 \index{Intermittency maps}
 \index{Neutral fixed points}
 Perhaps the simplest way to relax the uniform expansivity condition 
 is to allow some fixed (or periodic) point \( p \) to have a neutral 
 eigenvalue, e.g in the one-dimensional setting \( |Df(p)| =1 \), 
 while still requiring all other vectors in all directions over the 
 tangent spaces of all points to be strictly expanded by the action 
 of the derivative (though of course not uniformly since the 
expansion 
 must degenerate near the point \( p \)). Remarkably this can have 
 extremely dramatic consequences on the dynamics. 
 
There are some recent results for 
 higher-dimensional systems \cites{PolYur01a, Hu01, Gou03} but
a more complete picture is available 
the one-dimensional setting and thus we 
 concentrate on this case.  An initial motivation for 
 these kinds of examples arose from the concept of 
 \emph{intermittency}  in fluid dynamics.  A class of one-dimensional 
 maps expanding everywhere except at a fixed point was introduced by 
 Maneville and Pomeau in \cite{ManPom80} as a model of intermittency 
 since numerical studies showed that orbits tend to spend a long time 
 \emph{trapped}  in a neighbourhood of the fixed point with 
relatively short 
 \emph{bursts of chaotic activity} outside this neighbourhood. 
 Recent work  shows that indeed, these long periods of inactivity 
 near the fixed point are a key to slowing down the mixing process 
and 
 obtaining examples of systems with subexponential decay of 
 correlations. 
  
We shall consider interval maps \( f \) which are 
piecewise \( C^{2} \) 
     with a \( C^{1} \) extension to the boundaries of the \( C^{2} 
     \) domains and for which the derivative is 
     strictly greater than 1 everywhere except at a 
  fixed point \( p \) (which for simplicity we can assume lies at the 
origin) 
  where \( Df(p)=1 \). 
  For definiteness, let us suppose that 
 on a small neighbourhood of 0 the map takes the form 
  \[ 
  f(x) \approx x+x^{2}\phi(x)
  \]
  where \( \approx \) means that the terms on the two sides of the 
  expression as well as their first and second order derivatives 
converge as \( x\to 0 \). We assume moreover that
 \( \phi \) is \( C^{\infty} \) for \( x\neq 0 \); 
  the precise form of \( \phi  \) determines the precise degree of 
\emph{neutrality} 
  of the fixed point, 
and in particular  affects the second derivative \( D^{2}f \). It 
turn out that it 
plays a crucial role in determining the mixing properties and even 
the very existence of an absolutely continuous invariant measure. 
For the moment 
 we assume also a strong Markov property: 
 each domain of regularity of \( f \)
 is mapped bijectively to the whole interval. 
The following result shows that the situation can be drastically 
different 
from the uniformly expanding case.

\begin{theorem}\cite{Pia80}\label{noacip}
    If \( f \) is \( C^{2} \) at \( p \) (e.g. \( \phi(x) \equiv 1 \))
    then \( f \) does 
    not admit any acip.
 \end{theorem}

Note that 
\( f \) has the same topological behaviour as a uniformly 
expanding map,  
 typical orbits continue to wander densely on the whole interval, but 
 the proportion of time which they spend in various regions tends to 
 concentrate on the fixed point, so that, asymptotically, typical 
 orbits spend all their time near \( 0 \). 
 It turns out that in this situation there 
 exists an infinite (\( \sigma \)-finite) 
 absolutely continuous invariant measure which 
 gives finite mass to any 
 set not containing the fixed point  and infinite 
 mass to any neighbourhood of \( p  \)
 \cite{Tha83}. 
 
 The situation  changes if 
  we relax the condition that \( f \) be \( C^{2} \) at \( p \) and 
  allow the second derivative \( D^{2}f(x) \) to diverge to 
  infinity as \( x\to p \). This means that the derivative 
  increases 
   quickly as one moves away from 
   \( p \) and thus nearby points are repelled at a 
   faster rate. This is a very 
  subtle change but it makes all the difference. 
\index{Decay of correlations ! Subexponential}  
  \begin{theorem}\label{neutral}
       If \( \phi \) is  of the form 
  \( \phi(x) = x^{-\alpha}\)  for some \( \alpha\in (0,1)  \), then
\cites{Pia80, HuYou95} \( f \) admits an ergodic acip \( 
	  \mu \) and   
      \cites{Iso99, LivSauVai98, You99, PolYur01a, Sar02, Hu04, 
Gou04a}
\( \mu \) is mixing with decay of correlations 
      \[ 
      \mathcal C_{n} = \mathcal O(n^{1-\frac{1}{\alpha}}).
      \]
 If  \( 
       \phi(1/x) = \log x \log^{(2)}x\ldots\log^{(r-1)}x 
(\log^{(r)}x)^{1+\alpha} \)
       for some \( r\geq 1,\) \(\alpha\in (-1,\infty) \) 
       where \( \log^{(r)}=\log\log\ldots\log \) 
       repeated \( r \) times, then  \cite{Hol04}
       \( f \) admits a mixing acip \( \mu \) with 
       decay of correlations
       \[ 
       \mathcal C_{n}=\mathcal O (\log^{(r)}n)^{-\alpha}.
       \]
      \end{theorem}

Thus, the existence of an absolutely continuous invariant measure as 
in the 
uniformly expanding case 
 has been recovered, but the exponential rate of decay of correlation 
has not. 
We can think of 
 the indifferent fixed point as having the effect of 
  \emph{slowing down} 
  this process by trapping nearby points for disproportionately long 
time.   
 The estimates in \cites{Sar02, Hu04, Gou04a} include
    lower bounds as well as upper 
    bounds. The approach in \cite{You99} using Markov induced maps 
    applies also to non-Markov cases and 
    \cite{Hol04} can also be generalized to these cases.

The proofs of Theorem \ref{neutral} do not use directly the fact that 
\( 
f \) is nonuniformly expanding. Indeed the fact that \( f \) is 
nonuniformly expanding  does not follow 
automatically from the fact that the map is expanding away from the 
fixed point \( p \). However we can use the existence of the 
\emph{acip} to show that this condition is satisfied. Indeed by 
Birkhoff's Ergodic Theorem,  typical points spend a large proportion 
 of time near \( p \) but also  a positive proportion of 
 time in the remaining part of the space. More formally, 
 by a simple
 application of Birkhoff's Ergodic Theorem to the function \( \log 
|Df(x)| \), 
we have that, for \( \mu \)-almost every \( x \), 
 \[ 
 \lim_{n\to\infty}\frac{1}{n}\sum_{i=0}^{n-1}\log |Df(f^{i}(x))| \to 
\int 
 \log|Df| d\mu > 0.
 \]
  The fact that \( \int \log
 |Df|d\mu > 0 \) follows from the simple observation that \( \mu \) is
 absolutely continuous, finite,   and that \( \log |Df| > 0 \) 
 except at the neutral fixed point.

 \section{One-Dimensional Maps with Critical Points}
 \label{critical}
 \index{One-dimensional maps with critical points}
  We now consider another class of systems which can also exhibit
 various rates of decay of correlations, but where the mechanism for
 producing these different rates is significantly more subtle. 
 The most general set-up is that of a piecewise smooth 
 one-dimensional map \( f: I \to I \)  with 
 some finite set \( \mathcal C \) of \emph{critical/singular 
 points} at which  \( Df = 0 \)  
or \( Df = \pm \infty \) and/or at 
which \( f \) may be discontinuous.  There are at 
least two ways to quantify the ``uniformity: of the expansivity of \( 
f \) in ways 
that get reflected in different rates of decay of correlations: 
\begin{itemize}
    \item 
To consider the rate of growth of the derivatives along the orbits 
of the critical points; 
\item 
To consider the \emph{average} 
rate of growth of the derivative along \emph{typical }orbits. 
\end{itemize} 
In this 
section we will concentrate on the first, somewhat more concrete,
approach and describe the 
main results which have been obtained over the last 20/25 years. 
We shall focus specifically on the smooth case since this is where 
most 
results have been obtained. Some partial generalization to the 
piecewise smooth case can be found in \cite{DiaHolLuz04}.
The 
second approach is somewhat more abstract but also more general since 
it extends naturally to the higher dimensional context where 
critical points are not so well defined and/or cannot play such a 
fundamental role. The main results in this direction will be 
described 
in Section \ref{general} 
below in the framework of a general theory of non-uniformly expanding 
maps.

\subsection{Unimodal maps}
We  consider
 the class of \( C^{3} \) interval maps \( f: I\to I\) with some
 finite set \( \mathcal C \) of non-flat critical points.  
We recall that  \( c
 \) is a critical point if \( Df (c) = 0 \); the critical
 point is non-flat if there exists an \( 0< \ell<\infty \) called the
 \emph{order} of the critical point, such that \( |Df(x)| \approx
 |x-c|^{\ell-1} \) for \( x \) near \( c \);  \( f \) is
 unimodal if it has only one critical point, and multimodal if it has
 more than one. 
Several results to be mentioned below have been proved under 
a standard technical  \emph{negative Schwarzian derivative} condition 
which is a kind of convexity assumption on the derivative of \( f \), 
see \cite{MelStr88} for details. Recent results \cite{Koz00} indicate 
that this 
condition is  often superfluous and thus we will not mention it 
explicitly.

The first result on the statistical properties of such maps 
 goes back to Ulam and von Neumann \cite{UlaNeu47} who
 showed that the \emph{top} unimodal quadratic map, \( f(x) =
x^{2}-2\) has an \emph{acip}. 
Notice that this map is actually a Markov map but does not 
satisfy the bounded distortion condition due to the presence of the 
critical point. It is possible to construct an induced 
 Markov map for \( f \) which does satisfy this condition and gives 
 the result, but Ulam and von 
Neumann used a very direct approach, observing that \( f \) is 
\( C^{1} \) conjugate to a
piecewise-linear uniformly expanding Markov map for which Lebesgue 
measure is invariant and ergodic. 
This implies that the pull-back of the Lebesgue measure by the 
conjugacy is an \emph{acip} for \( f \). 
However, the existence of a smooth conjugacy is extremely rare and 
such 
an approach is not particularly effective in general. 

More general and more powerful approaches have allowed the existence 
of an \emph{acip} to be proved under increasingly general assumptions 
on the behaviour of the critical point. Let 
 \[
 D_{n}(c) = |Df^{n}(f(c))|. 
 \]
Notice that the derivative along the critical orbit needs to be 
calculated starting from the critical value and not from the critical 
point itself for otherwise it would be identically 0. 
 \index{Unimodal maps}
\begin{theorem}\label{uni}
    Let \( f: I \to I \) be a unimodal map. Then \( f \) admits an 
    ergodic acip if the following conditions hold (each condition 
    is implied by the preceding ones): 
    \begin{itemize}
	\item The critical point is pre-periodic \cite{Rue77};
\item 	The critical point is non-recurrent  \cite{Mis81};
\item \( D_{n} \to \infty \) exponentially fast \cites{ColEck83, 
NowStr88};
\item \( D_{n} \to \infty \) sufficiently fast so that \( \sum_{n} 
D_{n}^{1/\ell} < \infty \)\cite{NowStr91};
\item \( D_{n} \to \infty \) \cite{BruSheStr03}. 
\end{itemize}
If \( D_{n} \to \infty \) exponentially fast then some power of \( f 
\) 
is mixing and exhibits exponential decay of correlations 
\cite{KelNow92} (and \cite{You92} with additional bounded 
recurrence assumptions on the critical point). 
 \end{theorem}
 
 Notice that the condition of 
\cite{BruSheStr03} is extremely weak. In fact they show that it 
is sufficient for \( D_{n} \) to be eventually bounded below by some 
constant depending only on the order of the critical point. 
However even this condition is not optimal as 
there are examples of maps for which \( \liminf D_{n} = 0 \) but 
which still admit an ergodic \emph{acip}. 
It would be interesting to know whether an optimal condition is even 
theoretically possible: it is conceivable that  
a complete characterization of maps admitting 
\mbox{acip}'s in terms of the behaviour of the critical point is not 
possible because other subtleties come into play. 
 
\subsection{Multimodal maps}
Given the number of people and research 
papers in one-dimensional dynamics it is remarkable that until very recently 
there were essentially no results at all on the existence of 
\emph{acip}'s for multimodal maps. A significant breakthrough was  
achieved by implementing the strategy of constructing induced Markov 
maps and estimating the rate of decay of the tail. This strategy 
yields also estimates for various rates of decay in the 
unimodal case and extends very naturally to the multimodal case. 
\index{Decay of correlations ! Exponential}
\index{Decay of correlations ! Subexponential}
\index{Multimodal maps}
\begin{theorem}[\cite{BruLuzStr03}]
\label{decay}
Let \( f \) by a multimodal map with a finite set of critical points 
of order \( \ell \) and suppose that 
\begin{equation*}
\sum_n D_n^{-1/(2\ell-1)} < \infty \
\end{equation*}
for each critical point \( c \). Then there exists an ergodic acip \( 
\mu \) for \( f \). 
Moreover some power of \( f \) is mixing and the correlation function 
decays at the following rates:
\begin{description}
  \item [Polynomial case]
 If there exists \( C>0, \tau > 2\ell - 1 \) such that
 $$
D_n(c) \geq C n^{\tau}, \ 
 $$
for all
$c \in \mathcal C$ and $n \geq 1$,
then, for any \( \tilde \tau < \tfrac{\tau-1}{\ell-1} - 1 \), we have
\[
\mathcal C_n = \mathcal{O} (n^{-\tilde \tau})
\]
 \item [Exponential case]
 If there exist $C, \beta > 0$ such that
 $$
 D_n(c) \geq C e^{\beta n}
 $$
for all $c \in \mathcal C$ and $n\geq 1$, 
then there exist $\tilde \beta
> 0$ such
that
 \[
 \mathcal C_n = \mathcal O
(e^{-\tilde\beta n}).
 \]
   \end{description}
\end{theorem}

These results gives previously unknown
estimates for the decay of correlations even for unimodal maps in the 
quadratic family.  For example they imply that the so-called 
\emph{Fibonacci maps} \cite{LyuMil93} 
exhibit decay of correlation at rates which are 
faster than any polynomial. It seems likely that these estimates are 
essentially 
optimal although the argument only provides upper bounds.
The general framework of \cite{Lyn04} 
applies to these cases to provide estimates for the decay of 
correlation for observable which satisfy weaker than H\"older 
conditions on the modulus of continuity. 
The condition for the existence of an \emph{acip} have recently 
been weakened to the 
summability condition \( \sum D_{n}^{-1/\ell}<\infty \) and to
allow the possibility of critical points of different orders 
\cite{BruStr01}. Some improved technical expansion estimates have been 
also obtained in \cite{Ced04} which allow the results on the decay of 
correlations to apply to maps with critical points of different orders.

Based on these results, 
the conceptual picture of the causes of slow rates of decay of 
correlations appears much more similar to the 
case of maps with indifferent fixed points 
than would appear at first sight: we can think of the case 
in which the rate of growth of \( D_{n} \) is subexponential as a 
situation in which the critical orbit is \emph{neutral} 
or \emph{indifferent} and points which land close to the critical 
point tend to 
remain close to (``trapped'' by) its orbit  for a 
particularly long time. During this time orbits  are 
behaving ``non-generically'' and are not 
distributing themselves over the whole space as uniformly as they 
should.  Thus the mixing process is delayed and the rate of decay of 
correlations is correspondingly slower. 
When \( D_{n} \) grows exponentially, the critical orbit can be 
thought of (and indeed \emph{is}) a non-periodic  \emph{hyperbolic 
repelling} 
orbit and nearby points are pushed away exponentially fast. Thus 
there is no significant loss in the rate of mixing, and the decay of 
correlations is not significantly slowed down
notwithstanding the presence of a critical point.

\subsection{Benedicks-Carleson maps}
\label{indMark}
\index{Induced Markov maps ! in the quadratic family}
 We give here a sketch of the construction 
 of the induced Markov map for a class of unimodal maps. 
 We shall try to give a conceptually clear description of the main 
 steps and ingredients required in the construction. 
 The details of the argument are unfortunately particularly technical 
 and a lot of notation and calculations are carried out only to 
 formally verify statements which are intuitively obvious. It is very 
 difficult therefore to be at one and the same time conceptually 
 clear and technically honest. We shall therefore concentrate 
 here on the former approach and make some remarks about the 
 technical details which we omit or present in a simplified form.

 \index{Quadratic family}
 Let 
 \[  f_{a}(x) = x^{2} - a \]
 for \( x\in I=[-2, 2] \) and 
 \[ 
a\in\Omega_{\varepsilon}=[2-\varepsilon, 2]
 \]
 for some \( \varepsilon > 0 \) sufficiently small.  
The assumptions and the details of the proof require the 
 introduction of several additional constants, some intrinsic to the 
 maps under consideration and some auxiliary for the purposes of the 
 argument. In particular we suppose that there is a \( \lambda \in 
 (0, \log 2)  \) and constants 
 \[ 
 \lambda \gg \alpha \gg  \hat\delta \gg \delta > 0 
 \] 
where \( x \gg y \) means that \( y \)  must be sufficiently 
 small relative to \( x \). Finally, to simplify the notation 
 we also let \( \beta = \alpha/\lambda \). 
We restrict ourselves to parameter values \( a\in 
\Omega_{\varepsilon} \) which satisfy the 
\emph{Benedicks-Carleson} conditions: 
\begin{description}
    \item[Hyperbolicity]
    There exist \( C > 0 \) such that 
 \[ 
 D_{n}\geq C e^{\lambda n}\quad \forall \ n\geq 1;
 \]
 \item[Slow recurrence]
  \begin{equation*}
      |c_{n}|\geq e^{-\alpha n} \quad \forall \ n\geq 1.
  \end{equation*}
\end{description}
In section \ref{existence} on page \pageref{existence} 
we sketch a proof of the fact that these 
conditions are satisfied for a positive measure set of parameters in 
\( \Omega_{\varepsilon} \) (for \emph{any} \( \lambda \in (0, \log 2) 
\) and 
\emph{any }\( \alpha > 0 \)). They are therefore reasonably generic 
conditions. Assuming them here will allow us 
 to present in a compact form an almost complete proof. During the 
 discussion we shall make some comments about how the argument can be 
 modified to deal with slower rates of growth of \( D_{n} \) and 
 arbitrary recurrence patterns of the critical orbit. 

 We remark that the overall strategy as well as several details of 
the 
  construction in the two arguments (one proving that the 
  hyperbolicity and slow recurrence conditions occur with 
  positive probability and the other proving that they imply the 
  existence of an \textsf{acip}) are remarkably similar. This 
  suggests a deeper, yet to be fully understood and exploited, 
  relationship between the structure of dynamical space and that of 
  parameter space.

 We let 
\[
\Delta=(\delta, \delta) \subset (-\hat\delta, 
\hat\delta) = \hat\Delta 
\] 
denote \( \delta \) and \( \hat\delta \) neighbourhood of this 
critical point \( c \). 
The aim is to construct a Markov induced map 
\[ 
F: \Delta \to \Delta.
\]
We shall do this in three steps. We first define an induced map \( 
f^{p}: \Delta \to I \) which is essentially based on the time during 
which points in \( \Delta \) shadow the critical orbit. The shadowing 
time \( p \) is piecewise constant on a countable partition of \( 
\Delta \) but the images of partition 
elements can be arbitrarily small. Then we define an induced map \( 
f^{E}: \Delta \to I \) which is still not Markov but has
the property that the images of partition elements are uniformly 
large. 
Finally we define the Markov induced map \( F^{R}: \Delta \to \Delta 
\) as required.

\subsection{Expansion outside \protect\(\Delta\protect\)}
Before starting the construction of the induced maps, we state a 
lemma which gives some derivative expansion estimates outside the 
critical neighbourhood \( \Delta \). 
\begin{lemma}\label{l:expout1}
 There exists a constant \( C > 0 \) independent of \( \delta 
 \) such that for \( \varepsilon > 0 \) sufficiently small, all \( 
 a\in \Omega_{\varepsilon}, f=f_{a}, x\in I \) and \( n\geq 1 \) such 
 that \(   x, f(x), .., f^{n-1}(x) \notin \Delta   \) we have 
 \[ 
 |Df^{n}(x)| \geq  \delta e^{\lambda n}
 \]
 and if, moreover, \( f^n(x)\in\hat \Delta \) and/or \( x \in f(\hat 
\Delta) \)
 then 
 \[ 
 |Df^{n}(x)| \geq  C e^{\lambda n}
 \]
 \end{lemma}
 Notice that the constant \( C \) and the exponent \( \lambda \) do 
  not depend on \( \delta \) or \( \hat\delta \). 
  This allows us to choose \( \hat\delta \) and \( \delta \) small in 
  the following argument without worrying about this affecting the 
  expansivity estimates given here. 
In general of course both the constants \( C \) and \( \lambda \) 
depend on the size of this neighbourhood and it is an extremely 
useful feature of this 
particular range of parameter values that they do not. 
In the context of the quadratic family these estimates can be proved 
 directly using the smooth conjugacy of the top map \( f_{2} \) with 
 the piecewise affine tent map, see \cites{UlaNeu47, Luz00}. However 
 there are general theorems in one-dimensional dynamics to the effect 
 that one has uniform expansivity outside an arbitrary neighbourhood 
 of the critical point under extremely mild conditions \cite{Man85} 
 and this is sufficient to treat the general case in 
 \cite{BruLuzStr03}.   

\subsection{Shadowing the critical orbit}
We start by defining a partition of the critical neighbourhoods \( 
\Delta \) and \( \hat\Delta \).  
For any integer $r\geq 1$ let 
 \(I_{r} = [e^{-r},e^{-r+1}) \) and \( I_{-r}=(-e^{-r+1},-e^{-r}]\) 
 and,   for each \( r\geq r_{\hat\delta}+1 \), let 
 \( \hat I_{r}=I_{r-1}\cup I_{r}\cup I_{r+1} \).  
 We can suppose without loss of generality that \( r_{\delta} = \log 
 \delta^{-1}\) and \( r_{\hat\delta} = \log 
 \hat\delta^{-1}\) are integers
 Then
 \[ 
 \Delta=\{0\}\cup \bigcup_{|r|\geq r_{\delta}+1} 
I_{r}0,\quad\textrm{and}\quad
 \hat\Delta=\{0\}\cup \bigcup_{|r|\geq r_{\hat\delta}+1} I_{r}. 
 \]
This is one of the minor technical points of which we do not give a 
completely accurate description. Strictly speaking, the distortion 
estimates to be given below require a further subdivision of each \( 
I_{r} \) into \( r^{2} \) subintervals of equal length. This
does not affect significantly any of the other estimates. A similar 
partition is defined in Section \ref{omega n} on page \pageref{omega 
n} in somewhat more detail. We remark also that the need for the two 
neighbourhoods \( \Delta \) and \( \hat\Delta \) will not become 
apparent in the following sketch of the argument. We mention it 
however because it is  a 
crucial technical detail: the region \( \hat\Delta\setminus\Delta \) 
acts as a \emph{buffer zone} in which we can choose to apply the 
derivative estimates of Lemma \ref{l:expout1} or the shadowing 
argument of Lemma \ref{bindexp} according to which one is more 
convenient in a particular situation.

\index{Binding}
Now let 
\[
  p(r)=\max\{k:  |f^{j+1}(x) - f^{j+1}(c)|\leq  e^{-2\alpha j} \ 
\forall\ \ 
  x\in \hat I_{r}, \ \forall \  j <    k\}
  \] 
This definition was essentially first formulated in \cite{BenCar85} 
and \cite{BenCar91}.  The key characteristic is that it guarantees a 
bounded distortion property which in turn allows us to make several 
estimates based on information about the derivative growth along the 
critical orbit. Notice that the definition in terms of \( \alpha \) 
is based crucially on the fact that the critical orbit satisfies the 
slow recurrence condition. We mention below how this definition can 
be 
generalized. 
  
  \begin{lemma}\label{bindexp}
   For all points \( x\in \hat I_{r} \) and \( p=p(r) \) we have 
  \[
      |Df^{p+1}(x)|= |Df(x)| \cdot |Df^{p}(x_{0})| 
      \geq e^{(1-7\beta)r} 
   \]
  \end{lemma}
Recall that \( \beta=\alpha/\lambda \) can be chosen arbitrarily 
small. 
\begin{proof}
  First of all, using the bounded recurrence condition, the 
definition 
  of the binding period, and arguing as in the distortion estimates 
  for  the uniformly expanding maps above, it is not difficult to 
show 
  that 
  there is a constant 
  $\mathcal D_{1}$, depending on \( \alpha \) but independent of $r$ 
and $\delta$, 
  such that for all $x_0, 
  y_0 \in f(\hat I_{r})$ and $1\leq k \leq p$,
  \begin{equation}\label{binddist}
  \left| \frac{Df^{k}(x_0)}{Df^{k}(y_0)} \right| \leq \mathcal D_{1}.
  \end{equation}
  Using the definition of \( p \) this implies
  \[
  e^{-2\alpha (p-1)} \geq |x_{p-1}-c_{p-1}| \geq 
  \mathcal D^{-1}_1 |Df^{p-1}(c_{0})| \ |x_{0}-c_{0}|
   \geq   \mathcal D^{-1}_1 e^{\lambda (p-1)}e^{-2r}
  \]
  and thus 
  \(\mathcal{D}_1 e^{-2\alpha p}e^{2\alpha}\geq 
   e^{\lambda p}e^{-\lambda}e^{-2r}\).
  Rearranging gives 
  \begin{equation}\label{bindlen}
  p + 1 \leq 
  \frac{\log\mathcal{D}_1 +2\alpha+2\lambda+2 r + 
   }
  {\lambda+2\alpha} \leq \frac{3r}{\lambda}
  \end{equation}
  as long as we choose $\delta$ so that $r_\delta$ is sufficiently 
  large in comparison to the other constants, none of which depend on 
  \( \delta\).
  Moreover 
  \[
    \mathcal D 
  e^{-2r} |Df^p(x_{0})| 
  \geq \mathcal D |x_{0}-c_{0}|\ |Df^p(x_{0})| 
  \geq |x_{p}-c_{p}| \geq  e^{-2\alpha p}
  \]
  and therefore, using \eqref{bindlen}, we have 
\(
  |Df^{p}(x_{0})| \geq  
    \mathcal D^{-1} e^{2r}e^{-2\alpha p} 
  \geq  \mathcal D^{-1}   e^{(2-\frac{6\alpha }{\lambda}) r}.
\)
  Since \( x\in \hat I_{r}\) we have \( |Df(x)| = 2|x-c|\geq 
  2e^{-(r+2)}  \)
  and therefore
\(
  |Df^{p+1}(x)| = |Df^{p}(x_{0})| \ |Df(x)| 
 \geq e^{-2 } \mathcal D^{-1} 
  e^{(2-6\beta) r}e^{-r }
\geq 
   e^{-2 } \mathcal D^{-1} 
  e^{(1-6\beta) r}.
\)
This implies the result
  as long as we choose $r_\delta$ large enough.
  \end{proof}
Thus we have a first  induced map 
 \( F_{p}: \Delta \to I  \) given by \( F_{p}(x) = f^{p(x)}(x) \) 
where \( p(x) = p(r)  \) 
 for \( x\in I_{\pm r} \) which is uniformly expanding. 
  Indeed  notice that \( Df^{p(x)}(x) \to \infty \) as \( x\to c \). 
However there is no reason for which this map should satisfy the 
Markov property and indeed, an easy calculation shows that 
the images of the partition elements are \( \sim e^{-7\beta r} \to 
0\)  
and thus not even of uniform size.

The notion of shadowing can be generalized without any assumptions on 
the recurrence of the critical orbit in the following way, see  
\cite{BruLuzStr03}: let \( \{\gamma_{n}\} \) be a monotonically 
decreasing sequence  with \( 1 > \gamma_{n} > 0 \) 
and \( \sum \gamma_{n} < \infty \). Then for \( x\in\Delta \), let 
\[ 
p(x) := \max\{p: |f^{k}(x) - f^{k}(c)| \leq 
\gamma_{k} |f^{k}(c) - c| \ \ \forall \ k\leq p-1\}.
\]
A simple variation of the distortion calculation used above 
shows that the summability of \( 
\gamma_{n} \) implies that \eqref{binddist} holds with this 
definition 
also. Analogous bounds on \( D_{n} 
\)  will reflect the rate growth of the 
derivative along the critical orbit. If the growth of \( D_{n} \) is 
subexponential, the binding period will 
last much longer because the interval \( |f^{k}(x) - f^{k}(c)|  \) is 
growing at a slower rate.
The generality of 
the definition means that it is more natural to 
define a partition \( I_{p} \)  as the ``level sets'' of 
the function \( p(x) \).The drawback is that we have much less 
control over 
the precise size of these intervals and their distance from the 
critical point. Some estimates of the tail \( \{x> p\} \) can 
be obtained and it turns out these these are closely related to the 
rate of growth of \( D_{n} \) and to those of 
the return time function for the final induced Markov map. This is 
because the additional two steps, the escape time and the return time 
 occur exponentially fast. Thus \emph{the only bottleneck 
is the delay caused by the long shadowing of the 
critical orbit}.

\subsection{The escape partition}

Now let \( J\subset I \) be an arbitrary interval (which could also 
be 
\( \Delta \) itself). We want to
 construct a partition \( \mathcal P \) of \( J \) and a 
 stopping time \( E: J\to \mathbb N \) constant on elements of \( 
 \mathcal P \) with the property that for each \( \omega\in\mathcal P 
 \), \( f^{E(\omega)}(\omega) \approx \delta \). We think of \( 
 \delta \) as being our definition of \emph{large scale}; we call  
 \( E(\omega) \) the \emph{escape} time of \( \omega \), we call the 
 interval \( f^{E(\omega)}(\omega) \) and \emph{escape interval}, 
 and call 
 \( \mathcal P \) the \emph{escape time} partition of \( J \).

 The construction is 
 carried out inductively in the following way. 
 Let \( k\geq 1 \) and suppose that the intervals with 
 $E<k$ have already
 been defined. Let $\omega$ be a connected component of the 
complement of 
 the set $\{E<k\}\subset J.$ We consider the various cases 
 depending on the position
 of $f^k(\omega).$
 If $f^k(\omega)$ contains $\hat\Delta\cup I_{r_{\hat\delta}}\cup 
I_{-r_{\hat\delta}}$
 then we subdivide $\omega$ into three subintervals satisfying the 
required properties, and
 let $E=k$ on each of them.
 If $f^k(\omega)\cap\Delta=\emptyset$ we do nothing. 
 If $f^k(\omega)\cap\Delta\neq\emptyset$
 but $f^k(\omega)$ does not intersect more than two adjacent 
$I_{r}$'s then we say that
 \( k \) is an \emph{inessential return}
 and define the corresponding \emph{return depth} by 
 \( r= \max\{|r|: f^k(\omega)\cap
 I_{r}\neq \emptyset \}. \) 
 If \( f^{k}(\omega)\cap\Delta\neq \emptyset \) and
 \( f^{k}(\omega)\) intersects at least three elements of \( \mathcal
 I \), then we simply subdivide \( \omega \) 
 into subintervals $\omega_r$ in such a way that each $\omega_r$ 
satisfies
 \(
 I_{r}\subset f^k(\omega_r)\subset\hat{I}_r.
 \)
 For $r>r_{\delta}$ we say that $\omega_r$ has an essential return at 
time $k$ with an
 associated return depth $r$. For all other $r$, the
 $\omega_r$ are escape intervals, and for these intervals we set 
 $E=k$. Finally we consider one more important case. If \( 
f^{k(\omega)} 
 \) contains \( \Delta \) as well as (at least) the two adjacent 
 partition elements then we subdivide as described above except for 
 the fact that we keep together that portion of \( \omega \) which 
 maps exactly to \( \Delta \). We let this belong to the escape 
 partition \( \mathcal P \) but, as we shall see, we also let it 
 belong to the final partition associated to the full induced Markov 
 map. 

This purely combinatorial algorithm is designed to achieve two 
things, 
neither one of which follows immediatly from the construction:
 \begin{enumerate}
     \item Guarantee uniformly bounded distortion on each partition 
     element up to the escape time;
     \item Guarantee that almost every point eventually belongs to 
the 
     interior of a partition element \( \omega\in\mathcal P \).
 \end{enumerate}
 We shall not enter into the details of the distortion estimates here 
 but discuss the strategy for shomwing that \( \mathcal P \) is a 
 partition mod 0 of \( J \). Indeed this follows from a much stronger 
 estimate concerning the tail of the escape time function:
there exists a constant $\gamma> 0$ such that for any interval \( 
 J \) with \( |J|\geq \delta \)  we have
 \begin{equation}\label{basetail}
 |\{x\in J: E(x) \geq n \}| = 
\sum_{\substack{\omega'\in\mathcal{P}\\E(\omega')\geq n}} 
 |\omega'| \lesssim e^{-\gamma n} |J|.
 \end{equation}

The argument for proving \eqref{basetail} 
revolves fundamentally around the combinatorial 
information defined in the construction. More specifically,  for 
$\omega\in
 \mathcal{P}$ let $r_1,r_2,\dots,r_s$ denote the sequence of return 
depths 
 associated to essential return times occurring before $E(\omega')$, 
and
 let $\mathcal E(\omega)=\sum_{i=0}^{s}r_s$. Notice that this 
sequence may be empty 
 if $\omega$ escapes without intersecting $\Delta$, in this case we 
set
 $\mathcal E(\omega)=0$. 
 We now split the proof into three steps: 
\subsubsection*{1) Relation between escape time and return depths} 
 The first observation is that the escape time is 
 bound by a constant multiple of the sum of the return depths: 
 there exists a \( \kappa \) depending only on \( \lambda \) such that
     \begin{equation}\label{ret}
 E(\omega) \lesssim \kappa \mathcal E(\omega).
 \end{equation}
Notice that a constant  \( T_{0} \) should also be added to take care 
of the case in which \( \mathcal E(\omega)=0 \), corresponding to the 
situation in which 
 \( \omega \) has an escape the first time that iterates of \( \omega 
 \) intersect \( \Delta \). Since \( \omega \) is an escape, it has a 
 minimum size and the exponential growth outside \( \Delta \) gives a 
 uniform bound for the maximum number of iterates within which such a 
 return must occur. Since this constant is uniform it does not play a 
 significant role and we do not add it explicitly to simplify the 
 notation. 
 For the situation in which \( \mathcal E(\omega) > 0 \)
it is sufficient to show that each essential return with return 
depth \( r \) has the 
next essential return or escape within at most \( \kappa r \) 
iterations. Again this follow from the  observation that 
the derivative is growing exponentially on average 
during all these iterations: we have exponential growth outside \( 
\Delta \) and also exponential growth on average during each 
complete inessential binding period. This implie \eqref{ret}. 
From  \eqref{ret} we then have 
\begin{equation}\label{sum}
\sum_{\substack{\omega\in\mathcal{P}(\omega)\\E(\omega)\geq n}} 
\!\!\!\! |\omega|  \lesssim
\!\!\!\!  
\sum_{\substack{\omega\in\mathcal{P}(\omega)\\\mathcal E(\omega)\geq 
{n}/{\kappa}}} \!\!\!\! |\omega|. 
\end{equation}
Thus it is enough to estimate the right hand side of \eqref{sum} 
which is saying that there is an exponentially small probability of 
having a large total accumulated return depth before escaping, i.e. 
most intevals escape after relatively few and shallow return depths. 
The strategy is perfectly naive and consists of showing that the size 
of interval with a certain return depth \( \mathcal E \) is 
exponentially 
small in \( \mathcal E \) and that there cannot be too many so that 
their 
total sum is still exponentially small. 

\subsubsection*{2) Relation between size of \protect\( \omega \protect\) and return 
depth}
The size of each partition element can be estimated in terms of the 
essential return depths in a very coarse, i.e. non-sharp, way which 
is 
nevertheless sufficient for your purposes.  The argument relies on 
the 
following observation. Every return depth corresponds to a return 
which is 
followed by a binding period. During this binding period there is a 
certain overall growth of the derivative. During the remaining 
iterates there is also derivative growth, either from being outside 
\( 
\Delta \) or from the binding period associated to some inessential 
return. 
Therefore a simple 
application of the Mean Value Theorem gives 
 \begin{equation}\label{small}
 |\omega|\lesssim e^{- \frac{1}{2}\lambda \mathcal E(\omega)}. 
\end{equation}

\subsubsection*{3) The cardinality of the \protect\( \omega \protect\) with a certain 
return depth}
It therefore remains only to estimate the cardinality of of the set 
of 
elements \( \omega\) which can have the same value of \( \mathcal E 
\). 
To do this, notice first of all that we have a bounded multiplicity
of elements of $\mathcal{P}(\omega)$ which can share exactly the same 
sequence of return depths. More precisely
this corresponds to the number of escaping intervals which can arise 
at any given time from
the subdivision procedure described above, and is therefore less than 
$r_{\hat\delta}$.
Moroever, every return depth is bigger than $r_{\delta}$ and therefore
 for a given sequence $r_1\ldots, r_s$ we must have $s\leq \mathcal 
E/r_{\delta}$. 
Therefore letting $\eta=r^{-1}_{\delta},$ choosing $\delta$ 
sufficiently small 
the result follows from the following fact: 
Let $N_{k,s}$ denote the number of sequences $(t_1,\ldots,t_s)$, 
$t_i\geq 1$ for all $i$, $1\leq i \leq s$, such that 
$\sum_{i=1}^{s}t_i = k$. Then, 
for all $\hat\eta>0$, there exists $\eta>0$ such that for any 
integers $s,k$ with
$s<\eta k$ we have
\begin{equation}\label{ecount}
N_{k,s}\leq e^{\hat{\eta}k}.
\end{equation}
Indeed, applying
 \eqref{ecount} we get that the total number of possible sequences is
 \(
 N_{k}=\sum_{s=1}^{\eta \mathcal E} N_{\mathcal E,s}\leq \eta 
\mathcal E e^{\hat\eta \mathcal E}\leq e^{2\hat\eta \mathcal E}. 
 \)
 Taking into account the multiplicity of the number of elements 
sharing the same
 sequence we get the bound on this quantity as 
 \(\leq r_{\hat\delta}e^{2\hat\eta \mathcal E}\leq
 e^{3\hat\eta \mathcal E}.\) Multiplying this by  \eqref{small} and 
 substituting into \eqref{sum} gives the result.

To prove \eqref{ecount}, notice first of all that
 $N_{k,s}$ can be bounded above by 
 the number of ways to choose 
 $s$ balls from a row of $k+s$ balls, thus partitioning the remaining 
 $k$ balls into at most $s+1$ disjoint subsets. Notice also that this 
 expression is monotonically increasing in $s$, and therefore
 $$N_{k,s}\leq
 \begin{pmatrix}
  k+s\\ s
  \end{pmatrix}
  =\begin{pmatrix}
  k+s\\ k
  \end{pmatrix}
\leq \begin{pmatrix}
 (1+\eta)k \\ k
 \end{pmatrix}
 =\frac{[(1+\eta)k]!}{(\eta k)! k!}.
 $$
 Using Stirling's approximation formula $k!\in 
 [1,1+\frac{1}{4k}]\sqrt{2\pi k}k^{k}e^{-k}$, we have
\(
 N_{k,s} \leq\frac{[(1+\eta)k]^{(1+\eta)k}}{(\eta k)^{\eta k}k^k}
 \leq (1+\eta)^{(1+\eta)k}\eta^{-\eta k}
 \leq \exp\{(1+\eta)k\log(1+\eta) -\eta k\log\eta\}
 \leq \exp\{((1+\eta)\eta -\eta \log\eta)k\}.
\)
 Clearly $(1+\eta)\eta -\eta \log\eta\to 0$ as $\eta\to 0$. 
This completes the proof of \eqref{basetail}.

 \subsection{The return partition}
Finally we need to construct the full induced Markov map. 
To do this we simply start with \( \Delta \) and construct the escape 
partition \( \mathcal P \) of \( \Delta \). Notice that this is a 
refinement of the binding partition into intervals \( I_{r} \). 
Notice also that the 
definition of this escape partition allows as a special case 
the possibility that \( f^{E(\omega)}(\omega)=\Delta \). In this case 
of course \( \omega  \) satisfies exactly the required properties 
and we let it belong by definition to the partition \( \mathcal Q \) 
and define the return time of \( \omega \) as \( R(\omega) = 
E(\omega) \). 
Otherwise we consider each escape interval \( J = 
f^{E(\omega)}(\omega) \) 
and use it as a starting interval for constructing and escape 
partition and escape time function. Again some of the partition 
elements constructed in this way will actually have returns to \( 
\Delta \). These we define to belong to \( \mathcal Q \) and let 
their 
return time be the sum of the two escape times, i.e. the total number 
of iterations since they left \( \Delta \), so that \( 
f^{R(\omega)}(\omega) = \Delta \). 
For those that don't return to \( \Delta \) we repeat the procedure. 
We claim 
that almost every point of \( \Delta \) eventually belongs to an 
element which returns to \( \Delta \) in a good (Markov) way at some 
point and that the tail estimates for the return time function are 
not 
significantly affected, i.e. they are still exponential.

The final calculation to support this claim is based on the following 
fairly intuitive observation. Once an interval \( \omega \) has an 
escape, it has reached large scale and  therefore 
it will certainly cover \( \Delta \) 
after some uniformly bounded number of iterations. In particular it 
contains some subinterval \( \tilde\omega\subset\omega \) 
which has a return to \( \delta \) with at most this uniformly 
bounded number of iterates after the escape. Moreover, and 
crucially, the proportion of \( \tilde\omega \) in \( \omega \) 
is uniformly bounded below, i.e. 
there exists a constant \( \xi> 0 \) independent of \( \omega \) such 
that 
\begin{equation}\label{xi}
    |\tilde{\omega}|\geq \xi |\omega|.
\end{equation}
Using this fact we are now ready to estimate the tail of return 
times, 
 $|\{\omega\in\mathcal{Q}\ |\ R(\omega)>n \}|$.  
 The argument is again based on taking into account some 
combinatorial 
 information related to the itinerary of elements of the final 
 partition \( \mathcal Q \). In particular we shall keep track of the 
\emph{ number of escape times} which occur before time \( n \) for 
all 
 elements whose return is greater than \( n \). 
 First of all we let 
 \begin{equation}
 \mathcal{Q}^{(n)} = \{ \omega\in\mathcal{Q} | R(\omega) \geq n \}.
 \end{equation}
 Then, for each 
 $1\leq i\leq n$ we let 
 \begin{equation}
 \mathcal{Q}^{(n)}_{i} = \{ \omega\in\mathcal{Q}^{(n)} | 
 E_{i-1}(\omega)\leq n < E_{i}(\omega) \}
 \end{equation}
 be the set of partition elements in \( \mathcal Q^{(n)} \) who have 
 exactly \( i \) escapes. Amongst those we distinguish those with a 
 specific escape combinatorics. More precisely, for 
  $(t_1,\ldots,t_i)$ such that $t_j 
 \geq 1$ and $\sum t_j = n$,  let
 \begin{equation}
 \mathcal{Q}^{(n)}_{i}(t_1,...,t_i) = \left\{ 
\omega\in\mathcal{Q}^{(n)}_i 
 | \sum_{j=1}^{k} t_j = E_k(\omega), 1 \leq k \leq i-1 \right\}.
 \end{equation}
We then fix some small \( \eta > 0 \) to be determined below 
and write
 \begin{equation}\label{sum0}
 |\{ \omega\in\mathcal{Q} | R(\omega) \geq n \}| 
 = \sum_{i\leq n} |\mathcal{Q}^{(n)}_{i}| = \sum_{i\leq\eta n} 
 |\mathcal{Q}^{(n)}_{i}| + \sum_{\eta n < i\leq n} 
 |\mathcal{Q}^{(n)}_{i}|. 
 \end{equation}
By \eqref{xi} we have 
$|\mathcal{Q}^{(n)}_{i}| \lesssim
 (1-\xi)^{i}$, which gives
 \begin{equation}\label{sum1}
 \sum_{\eta n < i\leq n} |\mathcal{Q}^{(n)}_{i}| \lesssim 
 \sum_{\eta n < i\leq n} (1-\xi)^{i} \lesssim
 (1-\xi)^{\eta n} \approx e^{-\gamma_\xi n}
 \end{equation}
 for some $\gamma_\xi > 0$.
 Now let $\omega\subset\tilde{\omega}\in\mathcal{P}_{E_i}$ be one of 
the 
 non-returning parts of an interval $\tilde{\omega}$ that had its 
 $i$'th escape at time $E_i$. 
 Note that 
 
$f^{E_{i}}(\mathcal{P}_{E_{i+1}}|_{\omega})=\mathcal{P}(f^{E_i}(\omega))$. 
 
Therefore
 \begin{equation}
 \sum_{\substack{\omega'\subset\omega,\\ E_{i+1}(\omega')\geq 
 E_{i}+n}} |\omega'| \lesssim e^{-\gamma n}|\omega|,
 \end{equation}
 where $\gamma$ is as in \eqref{basetail}. 
 Let $\mathcal{Q}^{(n)}_{i}$ denote the set of intervals 
 $\omega\in\mathcal{Q}$ that have  precisely \( i \) escapes before 
 time \( n \) then
 \begin{equation}
 \sum_{\omega\in\mathcal{Q}^{(n)}(E_1,\ldots,E_i)} |\omega| \lesssim 
 e^{-\gamma n} |\Delta|.
 \end{equation}
 Therefore using again the combinatorial counting argument and the 
 inequality \eqref{ecount} we get 
 \begin{equation}\label{sum2}
 \begin{split}
 \sum_{i\leq\eta n} |\mathcal{Q}^{(n)}_{i}| &= \sum_{i\leq\eta n} 
 \sum_{(t_1,\ldots,t_i)} |\mathcal{Q}^{(n)}_{i}(t_1,...,t_i)|\\
 & \lesssim \sum_{i\leq\eta n} N_{n,i} e^{-\gamma n} |\Delta| 
 \lesssim e^{\hat{\eta}n} e^{-\gamma n}. 
 \end{split}
 \end{equation}
 Recall that by \eqref{ecount} \( \hat\eta \) can be chosen 
 arbitrarily small by choosing \( \eta \) small. Thus, 
 combining \eqref{sum1} and \eqref{sum2} and substituting into 
 \eqref{sum0} we get 
 \begin{equation*}
 |\{ \omega\in\mathcal{Q} | R(\omega) \geq n \}| \lesssim
 e^{-\gamma_\xi n} +  e^{(\hat{\eta} - 
 \gamma)n} \lesssim e^{-cn}. 
 \end{equation*}

\section{General Theory of Nonuniformly Expanding Maps}
\label{general}
\index{Nonuniformly expanding}

The intuitive picture which emerges from the examples discussed above 
is that of a \emph{default} exponential mixing rate for 
uniformly expanding systems  and H\"older continuous observables. 
However it is clear that general nonuniformly expanding systems 
can exhibit a variety of rates of decay. Sometimes these rates can be 
linked to properties of specific \emph{neutral} orbits 
which can \emph{slow down} the mixing process.   However 
it is natural to ask whether there is some intrinsic 
information related to the very definition of nonuniform expansivity
which determines the rate of decay of correlation. We recall 
that \( f \) is nonuniformly expanding if there
exists \( \lambda>0 \) such that for almost every \( x\in M \)
\begin{equation*} \tag*{\( (*) \)} 
\liminf_{n\to\infty}\frac{1}{n}
\sum_{i=0}^{n-1}\log \|Df^{-1}_{f^{i}(x)}\|^{-1}>\lambda.
\end{equation*}
Although the constant \( \lambda>0 \)  is uniform for Lebesgue almost 
every point, the convergence to the \( \liminf \)  is not generally
uniform. 

A \emph{measure of nonuniformity} has been proposed in 
\cite{AlvLuzPin03} based based precisely 
on the idea of quantifying the rate of convergence. 
 The measure has been shown to be directly 
linked to the rate of decay of correlations in \cite{AlvLuzPindim1} 
in the one-dimensional setting and in \cite{AlvLuzPin} in arbitrary 
dimensions, in the case of polynomial rates of decay. Recently the 
theory has been very extended to cover the exponential case as well 
\cite{Gou04}. We give here the precise statements. 

\subsection{Measuring the degree of nonuniformity}

\subsubsection{The critical set}
Let \( f: M\to M \)  be a (piecewise) \( C^{2} \) map. 
For \( x\in M \) we let \( Df_{x} \) denote the derivative of \( f \) 
at \( x \) and define
\(
\|Df_{x}\| = \max\{\|Df_{x}(v)\|: v\in T_{x}M, \|v\|=1\}.
\) 
We suppose that \( f \) fails to be 
a local diffeomorphism on some zero measure \emph{critical} set 
\( \mathcal C \) at which \( f \) may be discontinuous and/or 
\( Df \) may be discontinuous and/or singular and/or blow up to 
infinity.  Remarkably, all these cases can be treated in a unified 
way 
as \emph{problematic} points as will be seen below. In particular we 
can define a natural 
generalization of the non-degeneracy (non-flatness) condition 
for critical points of one-dimensional maps. 

\begin{definition}
 The \emph{critical set} \( \mathcal C \subset M \) is
 \emph{non-degenerate} if \( m(\mathcal C)=0 \) 
 and there is a constant $\beta>0$
 such that for every $x\in M\setminus\mathcal C$ we have
 $\dist(x,\mathcal C)^{\beta}\lesssim \|Df_{x}v\|/\|v\|\lesssim
 \dist(x,\mathcal C)^{-\beta}$ for all $v\in T_x M$, and the 
functions \(
 \log\det Df \) and \( \log \|Df^{-1}\| \) are \emph{locally 
Lipschitz}
 with Lipschitz constant \( \lesssim \dist(x,\mathcal C)^{-\beta}\).
 \end{definition}
 
From now on we shall always assume these non-degeneracy conditions. 
We remark that the results to be stated below are non-trivial 
even when the critical set \( \mathcal C \) is empty and 
\( f \) is a local diffeomorphism everywhere. 
For simplicity we suppose also that \( f \) is topologically 
transitive, i.e. there exists a point \( x \) whose orbit is dense in 
\( 
M \). 
Without the topological transitivity condition we would just get 
 that the measure \( \mu \) admits a finite number of ergodic 
  components and the results to be given below would then apply to 
  each of its components. 
    
  \subsubsection{Expansion and recurrence time functions}
  Since we have no geometrical information 
about \( f \) we want to show that the statistical 
    properties such as the rate of decay of 
    correlations somehow depends on abstract information related to 
    the non-uniform expansivity condition only. Thus we make the 
    following 
    \begin{definition}
    For \( x \in M \), we define the \emph{expansion time function} 
    \[
\mathcal E(x) =
\min\left\{N: \frac{1}{n}
\sum_{i=0}^{n-1} \log \|Df^{-1}_{f^{i}(x)}\|^{-1} \geq \lambda/2
\ \ \forall n\geq N\right\}.
\]
\end{definition}
By condition \( (*) \) this function is defined and finite almost 
everywhere. It measures the amount of time one has to wait before the 
uniform exponential growth of the derivative kicks in. If \( \mathcal 
E(x) \) was uniformly bounded, we would essentially be in the 
uniformly 
expanding case. In general it will take on arbitrarily large values 
and not be defined everywhere. 
If \( \mathcal E(x) \) is \emph{large} only on a \emph{small} set of 
points, then it makes sense to think of the map as being not very 
non-uniform, whereas, if it is large on a large set of points it is 
in some sense, very non-uniform. 
We remark that the choice of \( \lambda/2 \) in the definition 
   of the expansion time function \( \mathcal E(x) \) is fairly 
   arbitrary and does not  affect the asymptotic rate estimates. Any 
   positive number smaller than \( \lambda \) would yield the same 
   results. 
   
We also need to assume some dynamical conditions concerning 
the rate of recurrence of typical points near the critical set. 
We let \( d_{\delta}(x,\mathcal C) \) denote the \( \delta 
\)-\emph{truncated}
distance from \( x \) to \( \mathcal C \) defined as \(
d_{\delta}(x,\mathcal C) = d(x,\mathcal C) \) if \( d(x,\mathcal C)
\leq \delta\) and \( d_{\delta}(x,\mathcal C) =1 \) otherwise.
\begin{definition}
We say that \( f \) satisfies the property of 
{\em subexponential
recurrence} to the critical set if
for any $\epsilon>0$ there exists $\delta>0$ such that for Lebesgue
almost every $x\in M$
\begin{equation} \label{e.faraway1}\tag{\( ** \)}
    \limsup_{n\to+\infty}
\frac{1}{n} \sum_{j=0}^{n-1}-\log \dist_\delta(f^j(x),\mathcal C)
\le\epsilon.
\end{equation}
\end{definition}
We remark that although condition \( (**) \) might appear to be a 
very technical 
condition, it is actually quite natural and in fact \emph{almost} 
necessary. 
Indeed, suppose that an absolutely continuous invariant 
measure \( \mu \) did exist for \( f \). Then, a simple application 
of 
Birkhoff's Ergodic theorem implies that condition \( (**) \) is 
equivalent to the integrability condition 
\[ 
\int_{M} |\log \dist_\delta(x,\mathcal S)| d\mu < \infty
\]
which is simply saying that the invariant measure 
does not give too 
much weight to a neighbourhood of the discontinuity set.

Again, we want to differentiate between 
different degrees of recurrence in a similar way to the way we 
differentiated between different degrees of non-uniformity of the 
expansion. 
\begin{definition}
       For \( x \in M \), we define the \emph{recurrence time 
function} 
\[
\textstyle{
\mathcal R(x) = \min\left\{N\ge 1: \frac{1}{n} \sum_{j=0}^{n-1} -\log
\dist_\delta(f^j(x),\mathcal C) \leq 2\varepsilon, \forall n\geq 
N\right\}
}
\]
\end{definition}
Then, for a map satisfying both conditions \( (*) \) and \( (**) \) 
we 
let 
\[ 
\Gamma_{n} = \{x: \mathcal E(x) > n \text{ or
} \mathcal R(x) > n\}
\]
Notice that \( \mathcal E(x) \) and \( \mathcal R(x) \) are finite 
almost everywhere and thus \( \Gamma_{n}\to 0 \). It turns out that 
the rate of decay of \( |\Gamma_{n}| \) is closely related to the 
rate of decay of correlations. In the statement of the theorem we let 
\( \mathcal C_{n} \) denote the correlation function for H\"older 
continuous observables. 

\index{Decay of correlations ! Exponential}
\index{Decay of correlations ! Subexponential}

\begin{theorem}
      Let \( f:M\to M \) be a transitive \( C^{2} \) local
 diffeomorphism outside a non-degenerate critical set \( \mathcal C 
 \),
 satisfying conditions \( (*) \) and \( (**) \). Then 
 \begin{enumerate}
\item  \cite{AlvBonVia00} 
\( f \) admits an acip \( \mu \). Some power of \( f \) is mixing
\item 
\cites{AlvLuzPin03, AlvLuzPin, 
AlvLuzPindim1} 
Suppose that there exists \(
 \gamma>0 \) such that 
 \[
 |\Gamma_n| = \mathcal O(n^{-\gamma}).  
 \] 
Then  
    \[
    \mathcal C_{n} = \mathcal O(n^{-\gamma+1 }).
    \]
    \item
    \cite{Gou04} 
    Suppose that there exists \( \gamma>0 \) such that 
    \[
    |\Gamma_n| = \mathcal O(e^{-\gamma n}).  
    \] 
   Then  there exists \( \gamma'>0 \) such that 
       \[
       \mathcal C_{n} = \mathcal O(e^{-\gamma' n }).
       \]
 \end{enumerate}      
    \end{theorem}

\subsection{Viana maps}
 \index{Viana maps}
 \index{Alves-Viana maps}
A main application of the general results described above are a class 
of maps known as \emph{Viana} or \emph{Alves-Viana} maps. 
Viana maps were introduced in \cite{Via97} as an example of a class 
of 
higher dimensional systems which are 
strictly \emph{not} uniformly expanding but for which the non-uniform 
expansivity condition is satisfied and, most remarkably, is
\emph{persistent} under small \( C^{3} \) perturbations, which is not 
the case for any of the examples discussed above. These maps are 
defined as skew-products on a two dimensional cylinder of the form
 \(  f: \mathbb S^1\times\mathbb R
\to \mathbb S^1\times \mathbb R \)
\[ 
f (\theta, x) = 
( \kappa \theta, x^{2}+a + \varepsilon \sin 2\pi\theta)
 \]
where \( \varepsilon \) is assumed sufficiently small and \( a \) is 
chosen so that the one-dimensional quadratic map \( x\mapsto x^{2}+a 
\) 
for which the critical point lands after a finite number of iterates 
onto a hyperbolic repelling periodic orbit (and thus is a \emph{good} 
parameter value and satisfies the non-uniform expansivity conditions 
as mentioned above). The map \( \kappa\theta \) is 
 taken modulo \( 2\pi \), and the constant \( \kappa \) is a 
positive integer which was required to be \( \geq 16 \) in 
\cite{Via97} 
although it was later shown in \cite{BuzSesTsu03} that any integer \( 
\geq 2 \) 
will work. The \( \sin \) function in the skew product can also be 
replaced by more general Morse functions.

\begin{theorem}
    Viana maps
    \begin{itemize}
	\item \cite{Via97} satisfy \( (*) \) and \( (**) \). In particular 
	they are nonuniformly expanding;
\item 
\cites{Alv00,AlvVia02}
are  topologically mixing and
have a unique ergodic acip (with respect to two-dimensional Lebesgue 
measure);
\item \cite{AlvLuzPin}
have super-polynomial decay of correlations: for any \( 
    \gamma > 0 \) 
    we have
    \[ \mathcal C_{n}=  \mathcal O(n^{-\gamma});
    \]
    \item\cites{BalGou03, Gou04} have stretched exponential decay of 
    correlations: there exists \( \gamma>0 \) such that 
    \[ 
     \mathcal C_{n} =\mathcal O
    (e^{-\gamma \sqrt n}).
    \]
    \end{itemize}
\end{theorem}

\newcommand{\Int}{\operatorname{Int}}
\newcommand{\supp}{\operatorname{supp}}

\newcommand{\nc}{\newcommand}
\newcommand{\cal}[2]{\ensuremath{\mathcal{#1}^{(#2)}}}
\nc{\tcal}[2]{\ensuremath{\mathcal{\tilde{#1}}^{(#2)}}}
\nc{\hcal}[2]{\ensuremath{\mathcal{\hat{#1}}^{(#2)}}}

\section{Existence of Nonuniformly Expanding Maps}
\label{existence}
\index{Parameter exclusion}
\index{Nonuniform expansivity ! verifying}
An important point which we have not yet discussed is the fact that 
the verification of the nonuniform expansivity assumptions is a 
highly 
non-trivial problem. For example, the verification that Viana 
maps are nonuniformly expanding is one of the main results of 
\cite{Via97}. Only in some special cases can the required 
assumptions  be verified directly and easily. 
The definition of nonuniform expansivity is in terms of 
asymptotic properties of the map which are therefore 
intrinsically not checkable in any given finite number of steps. 
The same is true also  for the derivative growth assumptions 
on the critical orbits of one-dimensional maps as in Theorems 
\ref{uni} on page \pageref{uni} and \ref{decay} on page 
\pageref{decay}. A perfectly legitimate question is therefore 
whether these conditions actually do occur for any map at all. 
Moreover, recent results suggest  that this situation is at best 
extremely rare in the sense that the set of one-dimensional maps 
which have attracting periodic orbits, and in particular do not have 
an \emph{acip}, 
is \emph{open and dense} in the space of all one-dimensional maps 
\cites{GraSwi97, Lyu97, Koz03, She, KozSheStr03}. 
However, this topological point of view is only one way of defining 
``genericity'' and it turns out that for general one-parameter 
families of one-dimensional maps, the set of parameters 
for which an \emph{acip} does exist can have \emph{positive Lebesgue 
measure}  
(even though it may be topologically \emph{nowhere dense}). 

We give here a fairly complete sketch (!) of the argument 
in a special case, giving the complete description of the 
combinatorial construction and just brief overview of how the 
analytic 
estimates are obtained.  
For definiteness and simplicity we focus on the family 
\[ 
f_{a}(x) = x^{2}-a
\]
for \( x\in I= [-2, 2] \) and 
\[ 
a\in \Omega_{\varepsilon} = [2-\varepsilon, 2]
\]
for some \( \varepsilon > 0 \). 

\begin{theorem}[\cite{Jak81}]
\label{main theorem} 
For every \( \eta > 0 \) there exists an \( \varepsilon > 0 \) 
and a set \( \Omega^{*}\subset \Omega_{\varepsilon} \)
such that for all \( a\in \Omega^{*} \), \( f_{a} \) admits an ergodic
absolutely continuous invariant probability measure, and such that
\[
|\Omega^{*}|> (1-\eta) |\Omega|_{\varepsilon} > 0. 
\] 
\end{theorem}

There exists several generalizations of this result 
for families of smooth maps
\cites{Jak81, BenCar85, Ryc88, ThiTreYou94, MelStr93, Tsu93a, Tsu93} 
and even to families with completely degenerate (flat) 
critical points \cite{Thu99} and to piecewise 
smooth maps with critical points \cites{LuzVia00, LuzTuc99}. 
The arguments in the proofs are all fundamentally of a probabilistic 
nature and the conclusions depend on the fact 
that if \( f \) is nonuniformly expanding for a large  number 
of iterates \( n \) then it has a ``high probability'' of being 
nonuniformly expanding for \( n+1 \) iterates. 
Thus, by successively deleting those parameters which fail to be 
nonuniformly expanding up to some finite number of iterates has 
to delete smaller and smaller proportions. Therefore a 
positive proportion survives all exclusions.

In section \ref{combstruct} we give the formal inductive construction 
of the set 
\( \Omega^{*}\). In sections \ref{expout} and 
\ref{binding} we prove the two main technical 
lemmas which give expansion estimates for orbit starting respectively 
outside and inside some critical neighbourhood.  In section 
\ref{posexp} we prove the inductive step in the definition of \( 
\Omega^{*} \) and in section \ref{posmeas} we obtain the lower bound 
on the size of \( |\Omega^{*}| \). 

The proof involves several constants, some intrinsic to the family 
under consideration and some auxiliary for the purposes of the 
construction. The relationships between these constants and the order 
in which they are chose is quite subtle and also crucial to the 
argument. However this subtlety cannot easily be made explicit in 
such 
a sketch as we shall give here. We just mention therefore that there 
are essentially only two intrinsic constants: \( \lambda \) which is 
the expansivity exponent outside some (in fact any) critical 
neighbourhood, and \( \varepsilon \) which is the size of the 
parameter interval \( \Omega_{\varepsilon} \). \( \lambda \) can be 
chosen \emph{first} and is essentially arbitrary as long as \( 
\lambda \in (0, \log 2) \); \( \varepsilon \) needs to be chosen 
\emph{last} to guarantee that the auxiliary constants can be chosen 
sufficiently small. The main auxiliary constants are \( \lambda_{0} 
\) which can be chosen arbitrarily in \( (0, \lambda) \) and which 
gives the \emph{target} Lyapunov exponent of the critical orbit for 
good parameters, and  
\[
\lambda \gg \alpha \gg \hat\delta=\delta^{\iota} \gg \delta  > 0
\] 
which are chosen in the order given and sufficiently 
small with respect to the previous ones. During the proof we will 
introduce also some ``second order'' auxiliary constants which depend 
on these. Finally we shall use the constant \( C>0 \) to denote a 
generic constant whose specific value can in different formulae.

\subsection{The definition of \protect\( \Omega^{*} \protect\)}
\label{combstruct}
 
We let \(c_{0}=c_{0}(a)= f_{a}(0) \) denote the 
 \emph{critical value} 
 of \( f_{a} \) and for \( i\geq 0 \), 
 \( c_{i}=c_{i}(a) = f^{i} (c_{0})
 \). 
 For \( n\geq 0 \) and \(
 \omega\subseteq\Omega \) let 
 \(
 \omega_{n}=\{c_{n}(a); a\in \omega\}
 \subseteq I.
 \)
 Notice that for \( a=2 \) the critical value maps to a fixed point. 
 Therefore iterates of the critical point for parameter values 
 sufficiently close to \( 2 \) remain in an arbitrarily small 
 neighbourhood of this fixed point for an arbitrarily long time.
 In particular it is easy to see that all the inductive assumptions 
to be formulated below hold 
for all \( k\leq N \) where \( N \) can be taken arbitrarily 
 large if \( \varepsilon \) is small enough. This observation will 
play 
 an important role in the very last step of the proof.

\subsubsection{Inductive assumptions}
\label{indassu} 
Let \( \Omega^{(0)}=\Omega\) and \( 
\mathcal{P}^{(0)}=\{\Omega^{(0)}\} \)
denote
the trivial partition of
\( \Omega \). Given \( n\geq 1 \) suppose that for each \( k \leq n-1
\) there exists a set \( \Omega^{(k)}\subseteq\Omega \) satisfying the
following properties.

\begin{description}
\item[Combinatorics]
For the moment we describe the combinatorial structure 
as abstract data, the geometrical meaning of this data will become
clear in the next section.
There exists a partition \(\mathcal {P}^{(k)}\)
of \( \Omega^{(k)} \) into intervals such that each
\( \omega\in\mathcal P^{(k)} \) has an associated \emph{itinerary}
constituted by the following information
To each \( \omega\in\mathcal P^{(k)} \) is associated a 
sequence \(
0=\theta_{0}<\theta_{1}<\dots<\theta_{r}\leq k\),  \( r=r(\omega)\geq 
0 \)
of \emph{escape times}. Escape times are divided into three 
categories, i.e.
{\it substantial}, {\it essential}, and {\it inessential}.
Inessential escapes possess no combinatorial feature and 
are only relevant to 
the analytic bounded distortion argument to be 
developed later.
Substantial and essential escapes play a role in splitting
itineraries into  segments in the following sense.
Let \(
0=\eta_{0}<\eta_{1}<\dots<\eta_{s}\leq k\) , \(s=s(\omega)\geq 0 \)
be the maximal sequence of substantial and essential escape times. 
Between any of the two \( \eta_{i-1} \)
and \( \eta_{i} \) (and between \( \eta_{s} \) and \( k \)) there is a
sequence \( \eta_{i-1}<\nu_{1}<\dots<\nu_{t}<\eta_{i}\), \( 
t=t(\omega, i)\geq 0 \) of \emph{essential return times} (or
\emph{essential returns}) and between any two essential returns \(
\nu_{j-1} \) and \( \nu_{j} \) (and between \( \nu_{t} \) and \(
\eta_{i} \)) there is a sequence 
\(\nu_{j-1}<\mu_{1}<\dots<\mu_{u}<\nu_{j}\),
\( u=u(\omega, i, j)\geq
0\) of \emph{inessential return times} (or \emph{inessential
returns}).  Following essential and inessential return (resp. escape)
there is a time interval \( [\nu_{j}+1, \nu_{j}+ p_{j}]
\ \text{ (resp.  \(
[\mu_{j}+1, \mu_{j}+ p_{j}] \) ) } \) with \( p_{j} > 0 \) called the
\emph{ binding period}. A binding period cannot
contain any return and escape times.  
Finally, associated to each
essential and inessential return time (resp. escape) is a positive 
integer \( r \)
called the \emph{return depth} (resp. \emph{escape depth}).  

\item[Bounded Recurrence]
We define the function 
\(
\mathcal {E}^{(k)}: \Omega^{(k)}\to \mathbb N 
\) 
which associates to each \(
a\in\Omega^{(k)} \) the total sum of all essential return depths of
the element \( \omega\in\mathcal {P}^{(k)} \) containing \( a \)
in its itinerary up to and including time \( k \). Notice that 
\( \mathcal {E}^{(k)} \) is constant on elements of \( 
\mathcal{P}^{(k)} \) by construction. Then, for all
\( a\in \Omega^{(k)} \) 
\begin{equation*}
\tag*{\( (BR)_{k} \)} 
\mathcal {E}
^{(k)}(a) \leq \alpha k 
\end{equation*}

\item[Slow Recurrence]
For all \( a\in \Omega^{(k)} \) and all \( i\leq k \) we have
\begin{equation*}
\tag*{\( (SR)_{k} \)} 
|c_i(a)| \geq e^{-\alpha i}
\end{equation*}
Notice that \( \alpha \) can be chosen arbitrarily small as long as 
\( \varepsilon \) is small in order for this to hold for all \( i\leq 
N \). 

\item[Hyperbolicity]
For all 
\( a\in \Omega^{(k)} \) 
\begin{equation}\tag*{\( (EG)_{k} \)}
|(f_{a}^{k+1})'(c_{0})|\geq C e^{\lambda_{0} (k+1)}.  
\end{equation}

\item[Bounded Distortion]
Critical orbits with the same
combinatorics satisfy uniformly comparable derivative estimates:
For every  \( \omega\in\mathcal  {P}^{(k)} \), every pair of 
parameter values \( a,b\in\omega \) and every \( j\leq \nu+p+1 \) 
where \( \nu \) is
the last return or escape before or equal to time \( k \) and \( p
\) is the associated binding period, we have 
\begin{equation}\tag*{\( (BD)_{k} \)}
\frac{|(f_{a}^{j})'(c_{0})|}{|(f_{b}^{j})'(c_{0})|} \leq \mathcal D
\quad \text{ and } \quad \frac{|c'_{j}(a)|}{|c_{j}'(b)|}\leq \mathcal
D
\end{equation}
Moreover if $k$ is a 
substantial escape 
a similar distortion estimate holds for all $ j \leq l$ ($l$ is the 
next  chopping time) replacing $\mathcal D$ by $\tilde{\mathcal D}$ 
and
\( \omega \) by any subinterval \( \omega'\subseteq\omega \) which
satisfies \( \omega'_l \subseteq\Delta^{+} \). 
In particular  for \( j\leq k \), the map \( c_{j}: \omega \to 
\omega_{j} = \{c_{j}(a) : a\in \omega \} \) is a bijection.  

\end{description}

\subsubsection{Definition of \protect\( \Omega^{(n)} \text{ and } 
\mathcal P^{(n)}
\protect\)}
\label{omega n}  

For \( r\in\mathbb N \), let
\( I_{r}= [e^{-r}, e^{-r+1}),   I_{-r}=-I_{r}  \) and define
\[ 
\Delta^{+}=\{0\}\cup \bigcup_{|r|\geq r_{\delta^{+}}+1} I_{r}
\quad\text{ and }\quad \Delta=\{0\}\cup \bigcup_{|r|\geq r_{\delta}+1}
I_{r}. 
\] 
where \( r_{\delta}= \log\delta^{-1}, 
r_{\delta^{+}}= \iota\log\delta^{-1} \). We can suppose
without loss of generality that \( r_{\delta},
r_{\delta^{+}}\in\mathbb N \).  For technical reasons related to the 
distortion calculation we also need to subdivide each \( I_{r} \)
into \( r^2\) subintervals of equal length. This defines partitions \(
\mathcal I, \mathcal I^{+} \) of \( \Delta^{+}\) with \( \mathcal I=
\mathcal I^{+}|_{\Delta} \). An interval belonging to either one of
these partitions is of the form \( I_{r, m} \) with \( m\in [1, r^2]
\). Let \( I_{r, m}^{\ell} \) and \( I_{r,m}^{\rho} \) denote the
elements of \( \mathcal I^{+} \) adjacent to \( I_{r,m} \) and let \(
\hat I_{r,m}= I^{\ell}_{r,m}\cup I_{r,m}\cup I_{r,m}^{\rho} \). If \(
I_{r,m} \) happens to be one of the extreme subintervals of \(
\mathcal I^{+} \) then let \( I_{r,m}^{\ell} \) or \( I_{r,m}^{\rho}
\), depending on whether \( I_{r,m} \) is a left or right extreme,
denote the intervals 
\( 
\left(-\delta^{\iota}-\frac{\delta^{\iota}}{(\log\delta^{-\iota})^{2}}, 
-\delta^{\iota}\right] \) or \(
\left[\delta^{\iota}, \delta^{\iota}+
\frac{\delta^{\iota}}{(\log\delta^{-\iota})^{2}}\right) \) 
respectively. 
We now use this partition to define a
refinement \(\hat{\mathcal {P}}^{(n)}\) of \( \mathcal {P}^{(n-1)} 
\). 
Let \( \omega\in\mathcal
P^{(n-1)} \).  We distinguish two different
cases. 

\begin{description}
    \item[Non-chopping times]
    We
say that \( n \) is a non-chopping time for \( \omega\in\mathcal
P^{(n-1)}  \) if one (or
more) of the following situations occur: 
(1) \(
\omega_{n}\cap \Delta^{+}=\emptyset \); 
(2) \( n \) belongs to the
binding period associated to some return or escape time \( \nu<n \) of
\( \omega
\); 
(3)  \( \omega_{n}\cap\Delta^{+}\neq \emptyset \) but \(
\omega_{n} \) does not intersect more than two elements of the
partition \( \mathcal I^{+} \). 
In all three cases we
let \( \omega\in \hat{\mathcal {P}}^{(n)} \). 
In cases (1) and (2) no additional combinatorial information is added 
to the itinerary of \( \omega \). In case (3), 
 if \( \omega_{n}\cap (\Delta\cup I_{\pm
r_{\delta}})\neq \emptyset \) (resp. 
$\omega_n\subset\Delta^+\setminus 
(\Delta\cup I_{\pm r_{\delta}}$), we say that \( n \) is
an inessential return time (resp. inessential escape time ) for 
\( \omega\in \hat{\mathcal {P}}^{(n)}\). We define
the corresponding depth by 
\( r= \max\{|r|: \omega_{n}\cap
I_{r}\neq \emptyset \}. \)

\item[Chopping times]
In all
remaining cases, i.e. if \( \omega_{n}\cap\Delta^{+}\neq \emptyset \)
and \( \omega_{n} \) intersects at least three elements of \( \mathcal
I^{+} \), we say that \( n \) is a chopping time for \( \omega\in 
\mathcal
{P}^{(n-1)} \).  We define a natural subdivision \begin{equation*} 
\omega =
\omega^{\ell}\cup\bigcup_{(r, m)}\omega^{(r, m)}\cup \omega^{\rho}.
\end{equation*} so that each \( \omega_{n}^{(r, m)} \) fully contains
a unique element of \( \mathcal I_{+} \) (though possibly extending to
intersect adjacent elements) and \( \omega_{n}^{\ell} \) and \(
\omega_{n}^{\rho} \) are components of \( \omega_{n}\setminus
(\Delta^{+}\cap \omega_{n}) \) with \( |\omega_{n}^{\ell}|\geq
\delta^{\iota}/(\log\delta^{-\iota})^{2} \) and 
\( |\omega_{n}^{\rho}| \geq \delta^{\iota}/(\log\delta^{-\iota})^{2} 
\).
If the connected components of \( \omega_{n}\setminus (\Delta^{+}\cup
\omega_{n}) \) fail to satisfy the above condition on their length we
just glue them to the adjacent interval of the form \( w_{n}^{(r, m)}
\).  By definition we let each of the resulting subintervals of \(
\omega \) be elements of \( \hat{\mathcal {P}}^{(n)} \). 
The intervals \(\omega^{\ell}, \omega^{\rho} \) and \( 
\omega^{(r, m)} \) 
with \( |r|<r_{\delta} \) are called \emph{escape components} and are 
said to
have an \emph{substantial escape} and \emph{essential escape} 
respectively
at time \( n \). The corresponding values of \( |r|<r_{\delta} \) are 
the associated \emph{essential escape depths}. All other intervals 
are 
said to have an \emph{essential return} at time \( n \) and the
corresponding values of \( |r| \) are the associated \emph{essential
return depths.}
We remark that partition
elements \( I_{\pm r_{\delta}} \) do not  belong to \( \Delta \) but
we still say that the associated intervals  \( \omega^{(\pm
r_{\delta}, m)} \) have a return rather than an
escape. 
\end{description}

This completes the definition of the partition \(\hat{\mathcal 
{P}}^{(n)}\) of \(
\Omega^{(n-1)} \) and of  the function \(\mathcal {E}^{(n)} \) on \( 
\Omega^{(n-1)} \). We 
define 
\begin{equation}\label{exc}
\Omega^{(n)}= \{a\in\Omega^{(n-1)}: \mathcal {E}^{(n)}
(a) \leq {\alpha}n \} 
\end{equation}
Notice  that \( \mathcal E^{(n)} \) is constant on elements of
\( \hat{\mathcal {P}}^{(n)}\). Thus \( \Omega^{(n)} \) is the union 
of 
elements of \( \hat{\mathcal P}^{(n)} \) and we can define
\[
\mathcal {P}^{(n)}=\hat{\mathcal P}^{(n)}|_{\Omega^{(n)}}.
\]
Notice that the combinatorics and the recurrence condition \( (BR)_{n}
\) are satisfied for every  $a\in\Omega^{(n)}$
by construction. 
In Section \ref{posexp} we shall prove
that conditions \( (EG)_{n}, (SR)_n, (BD)_{n} \)
all hold for \( \Omega^{(n)} \). 
Then we define 
\[ 
\Omega^{*}=\bigcap_{n\geq 0}\Omega^{(n)}.  
\]
In particular, for every \( a\in\Omega^{*} \), the map \( f_{a} \) 
has 
an exponentially growing derivative along the critical orbit and 
thus, 
in particular, by Lemma \ref{uni} admits an ergodic \emph{acip}.  
In Section \ref{posmeas} we prove 
that  \( |\Omega^{*}|>0 \). 

We recall that a sketch of the proof of the existence of an induced 
Marov map under precisely the hyperbolicity and slow recurrence 
assumptions given here is carried out in section \ref{indMark}. 
As mentioned there, the strategy for construction of the induced 
Markov map is 
remarkably similar to the strategy for the construction carried out 
here for estimating 
the probability that such conditions hold. The deeper meaning of this 
similarity is not clear.

\subsection{Expansion outside the critical neighbourhood}
\label{expout}
On some deep level, the statement in the Theorem depends essentially 
on the following result which we have already used in section 
\ref{indMark}. 

\begin{lemma}\label{l:expout}
 There exists a constant \( C > 0 \) independent of \( \delta 
 \) such that for \( \varepsilon > 0 \) sufficiently small, all \( 
 a\in \Omega_{\varepsilon}, f=f_{a}, x\in I \) and \( n\geq 1 \) such 
 that \(   x, f(x), .., f^{n-1}(x) \notin \Delta   \) we have 
 \[ 
 |Df^{n}(x)| \geq  \delta e^{\lambda n}
 \]
 and if, moreover, \( f^n(x)\in\Delta^+ \) and/or \( x \in 
f(\Delta^{+}) \)
 then 
 \[ 
 |Df^{n}(x)| \geq  C e^{\lambda n}
 \]
 \end{lemma}
In the proof of the theorem we will use some other features of the 
quadratic family and of the specific parameter interval \( 
\Omega_{\varepsilon} \) but it is arguable that they are inessential 
and that the statement of Lemma \ref{l:expout} are to a certain 
extent sufficient conditions for the argument. It would be very 
interesting to try to prove the main theorem using only the 
properties stated in Lemma \ref{l:expout}. On a general
``philosophical'' level, the idea, as I believe originally 
articulated explicitly by L.-S. Young, is 
that 

\medskip
\noindent
\begin{center}
   \emph{Uniform Hyperbolicity for all parameters in most of the 
state 
    space \\ 
    implies 
    \\
    Nonuniform Hyperbolicity in all the state space for 
    most parameters.}
\end{center}

 \subsection{The binding period}
 \label{binding}

Next we make precise the definition of the \emph{binding 
period} which is part of the combinatorial information given above. 
Let \( k\leq n-1 \), \( \omega\in\mathcal P^{(k)}\) and 
suppose that \( k \) is an essential 
or inessential return or escape time for \( \omega \) with return 
depth \( r \). 
Then we define 
the \emph{binding period} of \( \omega_{k} \) as 
\[ 
p(\omega_{k}) = \min_{a\in \omega}\{p(c_{k}(a))\}
\]
where 
\begin{equation}\label{bindcond}
p(c_{k}(a)) = \min\{i: |c_{k+1+i}(a) - c_{i}(a)| \geq e^{-
2\alpha i}\}.
\end{equation}
This is the time for which the future orbit of \( c_{k}(a) \) can be 
thought of as \emph{shadowing} or being \emph{bound to} the orbit of 
the 
critical point (that is, in some sense, the number of iterations for 
which the orbit of \( c(a) \) repeats its early history after the \( 
k 
\)'th iterate). We will obtain some estimates concerning the length 
of this 
binding period and the overall derivative growth during this time.

\begin{lemma} \label{uniformbind} 
    There exist constants \( \tau_{0} > 0 \) 
   and \( \gamma_{1}\in (0,1) \) such that the following holds. 
   Let \( k\leq n-1 \), \( \omega\in\mathcal P^{(k)}\) and 
    suppose that \( k \) is an essential 
    or inessential return or escape time for \( \omega \) with return 
depth \( r \). 
  Let \( p=p(\omega_{k}) \). Then
for  every \( a\in\omega \)  we have 
\begin{equation}\label{up}
p \leq \tau_{0}\log |c_{k}(a)|^{-1} < k
\end{equation}
and
\begin{equation}\label{unibindder}
    |Df^{p+1}(c_{k}(a))| \geq   C e^{r (1-\gamma_{1})} 
    \geq C e^{\frac{1-\gamma_{1}}{\tau_{0}}(p+1)}
 \end{equation}   
and,  if \( k \) is an essential return  or an essential escape,  
then 
\begin{equation}\label{bindexp2}
|\omega_{k+p+1}|\geq  C e^{-\gamma_{1}r}.
\end{equation}
\end{lemma}
  
 To simplify 
the notation we write \( x=c_{k}(a) \) and \( x_{0}=c_{k+1}(a) \) 
and omit the dependence on the 
parameter \( a \) where there is no risk of confusion. 
The first step in the proof is to obtain a \emph{bounded distortion} 
estimate during binding periods:
there exists a constant \( \mathcal D_{1}(\alpha_{0}, \alpha_{1}) \) 
independent of \( x \), 
such that for all \( a\in \omega \), all \( y_{0}, z_{0}\in [x_{0}, 
c_{0}] \) 
and all \( 0\leq
j\leq \min\{p-1, k\} \) we have
\[
\left|\frac{(f^{j})'(z_{0})}{(f^{j})'(y_{0})}\right| 
< \mathcal
D_{1}.
\] 
This follows from the standard distortion calculations as 
in \eqref{eq:distest} on page \pageref{eq:distest}, using the 
upper bound \( e^{-2\alpha i} \) from the definition of binding
in the numerator and the lower bound \( e^{-\alpha i} \) from the 
bounded recurrence condition \( (SR)_{k} \) in the denominator. 
Notice that for this reason the
distortion bound is formally calculated for 
iterates \( j\leq \min\{p-1, k\} \) (the bounded recurrence condition 
cannot be guaranteed for iterates larger than \( k\)). 
The next step however gives an estimate for the duration of the 
binding period and implies that  \( p< k \) and therefore the 
distortion estimates do
indeed hold throughout the duration of the binding period. 
The basic idea for the upper bound on \( p \) is simple. 
The length of the interval \( [x_{0}, c_{0}] \) is determined 
by the length of the interval \( [x,c] \) which is \( c_{k}(a) \). 
The exponential growth of the derivative along the critical orbit and 
the bounded distortion imply that this interval is growing 
exponentially fast. The condition which determines the end of the 
binding period is shrinking exponentially fast. Some standard mean 
value theorem estimates using these two facts give the result. 
Finally the average derivative growth during the binding is given by 
the 
combined effect of the small derivative of order \( c_{k}(a) \) at 
the return 
to the critical neighbourhood and the exponential growth during the 
binding period. The result then intuitively boils down to showing 
that the 
binding period is long enough to (over) compensate the small 
derivative at the return. 

The final statement in the lemma requires some control over the 
way that the derivatives with respect to the parameter are related to 
the standard derivatives with respect to a point. This is a 
fairly important point which will be used again and therefore we give 
a more formal statement. 

\begin{lemma}\label{pardist} 
 There exists a constant \( \mathcal D_{2}>0 \) such that 
 for any \( 1\leq k \leq  n-1, \omega\in\mathcal P^{(k-1)} \) and \( 
a\in\omega \)
 we have
 \[ 
\mathcal D_{2}\geq \frac{|c_{k}'(a)|}{|Df^{k}_{a}(c_{0})|} \geq
\mathcal D_{2}
\]
and, for all \( 1\leq i < j\leq k+1 \), there exists \(
\tilde a\in\omega \) such that
\begin{equation}\label{pardisteq}  
\frac{1}{\mathcal D_{2} }
|Df_{\tilde a}^{j-i}
(c_{i}(\tilde a))| \leq 
\frac{|\omega_{j}|}{|\omega_{i}|}\leq 
\mathcal D_{2} |Df_{\tilde a}^{j-i}(c_{i}(\tilde a))|
\end{equation}
\end{lemma}
\begin{proof} 
 The second statement is a sort of \emph{parameter mean value 
 theorem} and follows immediately from the first one and the standard 
 mean value theorem. To prove the first one  let
\( F: \Omega \times I\to I \) be the function of two variables defined
inductively by \( F(a, x) = f_{a}(x) \) and \( F^{k}(a,
x)= F(a, f^{k-1}_{a}x).  \) Then, for \( x=c_{0} \), we have
\[
c_{k}'(a) = \partial_{a}F^{k}(a, c_{0}) = \partial_{a}F(a,
f_{a}^{k-1} c_{0}) = -1 + f'_{a}(c_{k-1}) c'_{k-1}(a).  
\]
Iterating
this expression gives
\begin{equation*}%
\begin{aligned}%
 - c_{k}'(a) = 1+
f_{a}'(c_{k-1})+f_{a}'(c_{k-1})f_{a}'(c_{k-2})+%
 \dots \\ \ldots %
+f_{a}'(c_{k-1})f_{a}'(c_{k-2}) \dots f'_{a}(c_{1})f'_{a}(c_{0})%
\end{aligned}%
\end{equation*} 
and dividing both sides by 
\(
(f^{k})'(c_{0})= f_{a}'(c_{k-1})f_{a}'(c_{k-2}) \dots
f'_{a}(c_{1})f'_{a}(c_{0}) \) gives %
\begin{equation}\label{parder}%
\frac{c_{k}'(a)}{(f^{k}_{a})'(c_{0})} = 
1+ \sum_{i=1}^{k}\frac{1}{(f^{i})'(c_{0})}.
\end{equation}
The result then depends on making sure that the sum on the right hand 
side is bounded away from -1. Since the critical point spends an 
arbitrarily large number \( N \) of iterates in an arbitrarily small 
neighbourhood of a fixed point at which the derivative is \( -4 \) we 
can bound an arbitrarily long initial part of this sum  by \( -1/2 
\). By the 
exponential growth condition the tail of the sum is still geometric 
and by taking \( N \) large enough we can make sure that this tail is 
less than \( 1/2 \) in absolute value. 
\end{proof}

Returning to the proof of Lemma \ref{uniformbind} we can use the 
parameter/space derivative bound to extend the derivative expansion 
result to the entire interval \( \omega_{k} \) and therefore to 
estimate the growth of this interval during the binding period.

\subsection{Positive exponents in dynamical space} 
\label{posexp}

Using a combination of the expansivity estimates outside 
\( \Delta \) and the binding period estimates for returns to \( 
\Delta 
\) it is possible to prove the inductive step stated above.

The slow recurrence condition is essentially an immediate 
consequence of the parameter exclusion condition. 

The exponential 
growth condition relies on the following crucial and non-trivial 
observation: \emph{the overall proportion of bound iterates is 
small}. This follows from the parameter exclusion condition which 
bounds the total sum of return depths (an estimate is required to 
show 
that inessential return do not contribute significantly to the total) 
and the binding period estimates which show that the length of the 
binding period is bounded by a fraction of the return depth. This 
implies that the overall derivative growth is essentially built up 
from the free iterates outside \( \Delta \) and this gives an overall 
derivative growth at an exponential rate independent of \( n \). 

The bounded distortion estimates again starts with the 
basic estimate as in \eqref{eq:distest} on page \pageref{eq:distest}.
By Lemma \ref{pardist} it is sufficient to prove the 
estimate for the space derivatives \( Df^{k} \); intuitively this is 
saying that critical orbits with the same
combinatorics satisfy the same derivative estimates. 
The difficulty here is that although the images of parameter 
intervals \( \omega \) are growing exponentially, they do not satisfy 
a uniform \emph{backward exponential} bound as required to carry out 
the step leading to \eqref{distest2} on page \pageref{distest2}; also 
images of \( \omega \) can come arbitrarily close to the critical 
point and thus the denominator does not admit any uniform bounds. 
The calculation therefore is technically quite involved and we refer 
the reader to published proofs such as \cite{Luz00} for the details. 
Here we just mention that the argument involves decomposing the sum 
into ``pieces'' corresponding to free and bound iterates and 
estimating each one independently, and taking advantage of 
the subdivision of the critical neighbourhood into interval \( I_{r} 
\) 
each of which is crucially further subdivided into further \( r^{2} 
\) subintervals of equal length. This implies that the contribution 
to 
the distortion of each return is at most of the order of \( 1/r^{2} 
\) 
instead of order 1 and allows us to obtain the desired conclusion 
using the fact that \( 1/r_{2} \) is summable in \( r \).

\subsection{Positive measure in parameter space}
\label{posmeas}

Recall that \hcal Pn is the partition of \( 
\Omega^{(n-1)} \) which takes into account the dynamics at time \( n 
\) 
and which restricts to the partition \( \mathcal P^{(n)} \) of \( 
\Omega^{(n)} \) after the exclusion of a certain elements of \( \hcal 
Pn \). Our aim here is to develop some combinatorial and metric 
estimates which will allow us to estimate the measure of parameters 
to 
be excluded at time \( n \). 

The first step is to take a fresh look at the combinatorial structure 
and ``re-formulate it'' in a way which is more appropriate. 
  To each \( \omega\in\hcal Pn \) is 
associated a sequence 
\( 0=\eta_{0}<\eta_{1}<\dots <\eta_{s}\leq n, \ s=s(\omega)\geq 0 
\) of escape times and 
a corresponding sequence of escaping components 
\(
	\omega\subseteq\omega^{(\eta_{s})}\subseteq\dots
	\subseteq\omega^{(\eta_{0})}
\)	\quad 
with 
\(\omega^{(\eta_{i})}\subseteq\Omega^{(\eta_{i})} \)
and 
\( \omega^{(\eta_{i})}\in\mathcal P^{(\eta_{i})}.\)
To simplify the formalism we also define some ``fake'' escapes by 
letting
\(
	\omega^{(\eta_{i})}=\omega
\) 
for all  \(s+1\leq i\leq n\). In this way we have a
well defined parameter 
interval \( \omega^{(\eta_{i})} \) associated to \( \omega\in
\hcal Pn \) for each \( 0\leq i\leq n \). 
Notice that for two intervals \( \omega, \tilde\omega\in \hcal Pn \) 
and 
any \( 0\leq i\leq n \), the corresponding intervals 
\( \omega^{(\eta_{i})} \) and \( \tilde\omega^{(\eta_{i})} \) 
are either disjoint or coincide.
Then we define
\[ 
Q^{(i)}=\bigcup_{\omega\in\hcal Pn}\omega^{(\eta_{i})}
 \] 
 and let 
 \(
 \cal Qi = \{\omega^{(\eta_{i})}\}
 \)
 denote the natural partition of 
 \( Q^{(i)} \) into intervals of the form \( \omega^{(\eta_{i})} \). 
 Notice that 
 \(
 \Omega^{(n-1)}= Q^{(n)}\subseteq\dots \subseteq
 Q^{(0)}=\Omega^{(0)} 
 \) 
 and 
  \(\cal Qn=\hcal Pn \)
 since
 the number \( s \) of escape times is always strictly 
 less than \( n \) and therefore in particular 
 \( \omega^{(\eta_{n})}=\omega \) 
 for all \( \omega\in\hcal Pn \).
 For a given \( \omega=\omega^{(\eta_{i})}\in \cal Qi, 
 \ 0\leq i \leq n-1 \) we let 
 \[ 
  Q^{(i+1)}(\omega) = \{\omega'=\omega^{(\eta_{i+1})}
 \in \cal Q{i+1}: \omega'\subseteq \omega\} 
 \] 
 denote all the elements of \( \cal Q{i+1} \) which are contained 
 in \( \omega \) 
 and let \( \cal Q{i+1}(\omega) \) denote the corresponding partition.
Then we define a function \( \Delta\cal Ei: Q^{(i+1)}(\omega) \to 
\mathbb N  \) by
\[ 
\Delta\cal Ei (a)=\cal E{\eta_{i+1}}(a)
-\cal E{\eta_{i}}(a). 
\]
 This gives the total sum of all essential return 
depths associated to the itinerary of the element \( \omega'
\in\cal Q{i+1}(\omega) \) containing \( a \), between the escape 
at time \( \eta_{i} \) and the escape at time \( \eta_{i+1} \).  
Clearly \( \Delta\cal Ei (a) \) is constant on elements of 
\( \cal Q{i+1}(\omega) \). 
Finally we let 
\[
 \cal Q{i+1}(\omega, R) = \{\omega'\in \mathcal Q^{(i+1)}: 
 \omega'\subseteq \omega, \Delta\cal Ei (\omega')= R\}.
 \]
 Notice that the entire construction given here depends on \( n \). 
The main motivation for this construction and is the following 

 \begin{lemma}
     There exists a constant \( \gamma_{0}\in (0, 1-\gamma_{1})  \) 
     such that the following holds. 
 For all \( i \leq n-1\), \( \omega\in \cal Qi \)  and \( R\geq 0 \) 
we have 
 \begin{equation}\label{main}
 \sum_{\tilde\omega\in \cal Q{i+1}(\omega, R)} \!\!\!\!\! 
|\tilde\omega|
 \leq e^{(\gamma_{1}+\gamma_{0} -1) R} |\omega|.
 \end{equation}
 \end{lemma}

 This says essentially that the probability of accumulating a large 
 total return depth between one escape and the next is exponentially 
 small. 
The strategy for proving this result is straightforward. We show 
first of all that 
 for  \(  0\leq i \leq n-1\), \( \omega\in \cal Qi \), \( R\geq 0 \) 
and 
 \(\tilde\omega\in \cal Q{i+1}(\omega, R)\) 
  we have
 \begin{equation}\label{metest}
 |\tilde\omega|\leq e^{(\gamma_{1}-1) R} |\omega|. 
 \end{equation}
The proof is not completely straightforward but depends on the 
intuitively obvious fact that an interval which has a deep return 
must necessarily be very small (since it is only allowed to contain 
at most three adjacent partition elements at the return). Notice 
moreover that this statement on its own is not sufficient to imply 
\eqref{main} as there could be many small intervals which together 
add up to a lot of intervals having large return. However we can 
control to some extent the multiplicity of these intervals 
and show that we can choose an arbitrarily small \( \gamma_{0} \) (by 
choosing the critical neighbourhood \( \Delta \) sufficiently small) 
so that 
    for all \( 0\leq i
\leq n-1\), 
\( \omega\in \mathcal{Q}^{(i)}\) and 
\( R\geq r_{\delta} \), 
we have 
\begin{equation} \label{combestlem}
\#
\mathcal{Q}^{(i+1)}_n(\omega, R) \leq e^{\gamma_{0} R}.  
\end{equation} 
This depends on the observation that each \( \omega \) has an 
essentially unique (uniformly bounded multiplicity) sequence of 
return 
depths. 
Thus the estimate can be approached via purely 
 combinatorial arguments very similar to those  used in 
relation to equation \eqref{ecount}. 
Choosing \( \delta \) small means the sequences 
of return depths have terms bounded below by \( r_{\delta} \) which 
can be chosen large, and this allows the exponential rate of increase 
of the combinatorially distinct sequences with \( R \) to be taken 
small. Combining \eqref{metest} and \eqref{combestlem} immediately 
gives \eqref{main}.

Now choose some \( \gamma_{2}\in (0, 1-\gamma_{0}-\gamma_{1}) \) and 
let 
\[
\gamma=\gamma_{0}+\gamma_{1}+\gamma_{2} > 0.
\]
For   \(  0\leq i \leq n-1\),
and \( \omega\in \cal Qi \), write
\[ 
\sum_{\omega'\in\cal Q{i+1}(\omega)} \hspace{-.5cm} 
e^{\gamma_{2}\Delta\cal Ei (\omega')} 
|\omega'|  =
\sum_{\omega'\in\cal Q{i+1}(\omega, 0)}\hspace{-.6cm}|\omega'| + 
\sum_{R\geq r_{\delta}} 
e^{\gamma_{2}R} \hspace{-.6cm}\sum_{\omega'\in\cal Q{i+1}(\omega, R)}
\hspace{-.6cm}|\omega'|. 
\]
By  \eqref{main} we then have 
\begin{equation}\label{i}
\sum_{\omega'\in\cal Q{i+1}(\omega, 0)}\hspace{-.5cm}|\omega'| + 
\sum_{R\geq r_{\delta}} 
e^{\gamma_{2}R} \sum_{\omega'\in\cal Q{i+1}(\omega, 
R)}\hspace{-.6cm}|\omega'| 
\leq
\left(1+  \sum_{R\geq r_{\delta}} \hspace{-.2cm}
 e^{(\gamma_{0}+\gamma_{1}+ \gamma_{2}-1) R}\right) |\omega| 
\end{equation}
Since \( \cal En = \Delta\cal E0 + \dots + \Delta\cal E{n-1} \) 
and \( \Delta\cal Ei \) is constant on elements of \( \cal Q{i} \) 
we can write  
\begin{align*} 
\sum_{\omega\in\cal Qn} e^{\gamma_{2}\cal En(\omega)}|\omega| 
=
&\sum_{\omega^{(1)}\in \cal Q1 (\omega^{(2)})} 
\hspace{-.8cm} e^{\gamma_{2}\Delta\cal E0 (\omega^{(1)})}
\hspace{-.8cm}
\sum_{\omega^{(2)}\in \cal Q{2}(\omega^{(3)})} 
\hspace{-.8cm} e^{\gamma_{2}\Delta\cal E1 (\omega^{(2)})}
\hspace{-.3cm}\dots 
\\ 
&\dots 
\hspace{-1cm} \sum_{\omega^{(n-1)}\in \cal Q{n-1}(\omega^{(n)})} 
\hspace{-1.2cm} e^{\gamma_{2}\Delta\cal E{n-1}(\omega^{(n-1)})}
\hspace{-.7cm}
\sum_{\omega=\omega^{(n)}\in \cal Q{n}} 
\hspace{-.7cm}  e^{\gamma_{2}\Delta\cal E{n-1}(\omega^{(n)})}|\omega|.
\end{align*}
Notice the \textit{nested} nature of the expression. Applying 
\eqref{i} 
repeatedly gives
\begin{equation}\label{average}
    \int_{\Omega^{(n-1)}}\hspace{-.5cm}e^{\gamma_{2}\cal En} = 
    \sum_{\omega\in\cal Qn} \hspace{-.2cm} e^{\gamma_{2}\cal 
En(\omega)} |\omega| 
\leq
\left(1+  \sum_{R\geq r_{\delta}} 
 e^{(\gamma-1) R}\right)^{n}|\Omega|. 
\end{equation}
The definition of \( \Omega^{(n)} \) gives 
\[ 
|\Omega^{(n-1)}|-|\Omega^{(n)}|=|\Omega^{(n-1)}\setminus\Omega^{(n)}|=
|\{\omega\in\hcal Pn=\cal Qn: e^{\gamma_{2}\cal En}\geq 
e^{\gamma_{2}\alpha n}\}| 
\]
and therefore using Chebyshev's inequality 
and \eqref{average} we have
\[
|\Omega^{(n-1)}|-|\Omega^{(n)}| 
\leq 
e^{-\gamma_{2} \alpha 
n}\int_{\Omega^{(n-1)}}\hspace{-.5cm}e^{\gamma_{2}\cal En}
\leq 
\left[e^{-\gamma_{2} \alpha}\left(1+  \sum_{R\geq r_{\delta}} 
 e^{(\gamma -1) R}\right)\right]^{n}
|\Omega|
\]
which implies
\begin{equation*} 
   |\Omega^{(n)}|\geq |\Omega^{(n-1)}|-\left[e^{-\gamma_{2} \alpha}
   \left(1+  \sum_{R\geq r_{\delta}} 
 e^{(\gamma-1) R}\right)\right]^{n} |\Omega|
\end{equation*}
and thus 
\[ 
|\Omega^{*}|\geq \left(1-\sum_{j=N}^{\infty}
\left[e^{-\gamma_{2} \alpha}\left(1+  \sum_{R\geq r_{\delta}} 
 e^{(\gamma -1) R}\right)\right]^{j} \right) |\Omega|.
\]
Choosing \( N \) sufficiently large, by taking \( \varepsilon \) 
sufficiently small, guarantees that the right hand 
side is positive.

\section{Conclusion}
    \label{final}
In this final section we make some concluding remarks and 
present some questions and open problems. 
    
 \subsection{What causes slow decay of correlations ?}
\index{Decay of correlations}
The general theory described in section \ref{general} 
is based on a certain way of quantifying the intrinsic nonuniformity 
of \( f \) which does not rely on identifying particular critical 
and/or neutral orbits. However, the conceptual picture according to 
which slow rates of decay are caused by a slowing down process 
due to the presence of neutral orbits can also be generalized. 
Indeed, the abstract formulation of the concept of a neutral orbit 
is naturally that of an orbit with a zero Lyapunov exponent. 
The definition of nonuniform expansivity implies that almost all 
orbits have uniformly positive Lyapunov exponents but this does not 
exclude the possibility of some other point having a zero Lyapunov 
exponent. It seems reasonable to imagine that a point with a zero 
Lyapunov exponent could \emph{slow down} the overall mixing process 
in 
a way which is completely analogous to the specific examples 
mentioned 
above. Therefore we present here, in a heuristic form, 
a natural conjecture. 

\begin{conjecture}
    Suppose \( f \) is non-uniformly expanding. 
    Then \( f \) 
    has exponential decay of correlations if and only  
    it has no orbits with zero Lyapunov exponent.
 \end{conjecture}
 
An attempt to state this conjecture in a precise way reveal several 
subtle points which need to be considered. We discuss some of these 
briefly. Let \( 
\mathcal M \) denote the space of all probability \( f \)-invariant 
measures \( \mu \) on \( M \) which satisfy the integrability 
condition 
\(
\int  \log \|Df_{x}\| d\mu < \infty.
\)
Then by standard theory, see also \cite{BarPes04}*{Section 5.8},
we can apply a version of Oseledet's Theorem for non-invertible maps 
which says 
that there exist constants \( \lambda_{1}, \ldots, \lambda_{k} \) 
with \( k\leq d \), and 
a measurable decomposition \( T_{x}M = E^{1}_{x}\oplus\dots\oplus 
E^{k}_{x} \)
of the tangent bundle over \( M \) such that the decomposition is 
invariant by the derivative and such that for all \( j=1,\ldots, k \) 
and for all non zero vectors \( v^{(j)}\in 
E^{j}_{x} \) we have 
\[ 
\lim_{n\to\infty}\frac{1}{n}\sum_{i=0}^{n-1}\log 
\|Df^{n}_{x}(v^{(j)})\| = 
\lambda_{j}.
\]
The constants \( \lambda_{1}, \ldots, \lambda_{k} \) are called the 
\emph{Lyapunov exponents} associated to the measure \( \mu \).  
The definition of nonuniform expansivity implies that all Lyapunov 
exponents associated to the \emph{acip} \( \mu \) are \( \geq \lambda 
\) 
and thus uniformly positive, but it certainly does not exclude the 
possibility that there exist some other (singular with respect to 
Lebesgue) invariant probability measure with some zero Lyapunov 
exponent. This is the case for example 
for the maps of Section 
\ref{almost expanding} for which 
 the Dirac measure on the indifferent fixed point has a zero Lyapunov 
 exponent. 

 Thus one way to state precisely the above conjecture is to claim 
 that \( f \) has exponential decay of correlations if and only if 
 all Lyapunov exponents associated to all invariant probability 
 measure in \( \mathcal M \) are uniformly positive. 
 Of course, a priori, there may also be some exceptional points, not 
 typical for any measure in \( \mathcal M \), along whose orbit the 
 derivative expands subexponentially and which therefore might 
 similarly have a slowing down effect. Also it may be that one zero 
 Lyapunov exponent along one specific direction may not have a 
 significant effect whereas a measure for which all Lyapunov 
 exponents were zero would. Positive results in the direction of this 
 conjecture include the remarkable observation that 
local diffeomorphisms  for which all Lyapunov exponents 
 for all measures are positive, must actually be uniformly expanding
 \cites{AlvAraSau03, Cao03, CaoLuzRioHyp} and thus in particular have 
 exponential decay of correlations. 
 Moreover, in the context of 
 one-dimensional smooth maps with critical points it is known that in 
 the unimodal case exponential growth of the derivative along the 
 critical orbit (the Collet-Eckmann condition) 
 implies uniform hyperbolicity on periodic orbits 
 \cite{Now88} which in turn implies that all Lyapunov exponents of 
all 
 measures as positive
 \cite{BruKel98} and the converse is also true \cite{NowSan98}.
Thus conjecture 1 is  true in the one-dimensional unimodal setting.   
 
We remark that the assumption of non-uniform expansion is 
crucial here. There are several examples of systems which 
have exponential decay of correlations but clearly have invariant 
measures with zero Lyapunov exponents, e.g. 
partially hyperbolic maps or maps obtained as time-1 maps of certain 
flows \cites{Dol98a, Dol98b, Dol00}. These examples however are 
not non-uniformly expanding, and are generally \emph{partially 
hyperbolic} 
which means that there are two continuous subbundles such that the 
derivative restricted to one subbundle has very good expanding 
properties or contracting properties and the other subbundle has the 
zero Lyapunov exponents. For reasons which are not at all clear, this 
might be \emph{better} from the point of view of decay of 
correlations 
than a situation in which all the Lyapunov exponents of the 
absolutely 
continuous measure are positive but there is some \emph{embedded} 
singular measure with zero Lyapunov exponent slowing down the mixing 
process. Certainly there is still a lot to be understood on this 
topic.

\subsection{Stability}
The results on the existence of nonuniformly expanding maps for open 
sets or 
positive measure sets of parameters are partly \emph{stability} 
results. 
They say that certain properties of a system, e.g. being nonuniformly 
expanding, are stable in a certain sense. We mention here two other 
forms of stability which can be investigated. 

\subsubsection{Topological rigidity}
\index{Topological rigidity ! in one-dimensional dynamics}
 The notion of (nonuniform) expansivity is, \emph{a priori}, 
completely 
 metrical: it depends on the differentiable structure of \( f \) 
 and most constructions and estimates related to nonuniform 
expansivity 
 require delicate metric distortion bounds.  However the statistical 
 properties we deduce (the existence of an \emph{acip}, the rate of 
 decay of correlation) are objects and quantities 
 which make sense in a much more 
 general setting. A natural question therefore is whether the metric 
 properties are really necessary or just very useful conditions and 
to 
 what extent the statistical properties might depend only on the 
 underlying topological structure of \( f \). We recall that two maps 
\( 
 f: M \to M \) and \( g: N \to N \) are \emph{topologically 
 conjugate}, \( f \sim g \), 
 if there exists a homeomorphism \( h: M \to N \) such that \( h\circ 
 f = g \circ h \). We say that a property of \( f \) is 
 \emph{topological} or depends only on the \emph{topological 
structure} 
 of \( f \) if it holds for all maps in the topological conjugacy 
 class of \( f \). 
 
 The existence of an absolutely continuous invariant measure is 
 clearly not a topological invariant in general: if \( \mu_{f} \) is 
 an \emph{acip} for \( f \) then we can define \( 
\mu_{g}=h^{*}\mu_{f} 
 \) by \( \mu_{g}(A) = \mu_{f}(h^{-1}(A)) \) which
gives  an invariant probability measure but not absolutely continuous 
 unless the conjugating homeomorphism \( h \) is itself absolutely 
 continuous. For example the map of Theorem 
 \ref{noacip} has no \emph{acip} even though it is topologically 
 conjugate to a \emph{uniformly} expanding Markov map. However it 
 turns out, quite remarkably, that there are many situations in the 
 setting of one-dimensional maps with critical points in which the 
 existence of an \emph{acip} is indeed a topological property 
 (although there are also examples in which it is not \cite{Bru98a}). 
 Topological conditions which imply the existence of an \emph{acip} 
 for unimodal maps were given in \cites{Bru94, San95, Bru98b}. In 
 \cite{NowPrz98} (bringing together results of \cites{NowSan98, 
 PrzRoh98}) it was shown that the exponential growth condition 
 along the critical orbit for unimodal maps (which in particular 
 implies the existence of an \emph{acip}, see Theorem \ref{uni}) 
 is a topological property. A counterexample to this result in the 
 multimodal case was obtained in \cite{PrzRivSmi03}. However it was 
 shown in \cite{LuzWan03} that in the general multimodal case, if all 
 critical points are \emph{generic} with respect to the \emph{acip}, 
 then the existence of an \emph{acip} still holds for all maps in the 
 same conjugacy class (although not necessarily the genericity of the 
 critical points). 
 
 We emphasize that all these results do not rely on showing that all 
 conjugacies in question are absolutely continuous. Rather they depend
 on the 
 existence of some topological property which forces the existence 
 of an \emph{acip} in each map in the conjugacy class. These 
 \emph{acip}'s are generally \emph{not} mapped to each other by the 
 conjugacy. 
 
 \subsubsection{Stochastic stability}
 \index{Stochastic stability}
 Stochastic stability is one way to formalize the idea that the 
 statistical properties of a dynamical systems are 
 \emph{stable} under small random perturbations. There are  
several positive results on stochastic stability for uniformly 
expanding \cites{You86b, BalYou93,Cow00} and nonuniformly expanding  
maps in dimension 1 \cites{BalVia96,AraLuzVia}
and higher \cites{Ara01,AlvVia02}
See \cite{Alv03} for a comprehensive treatment of the results. 

\subsection{Nonuniform hyperbolicity and induced Markov maps}
The definition of nonuniform hyperbolicity in terms of conditions 
\( (*) \) and \( (**) \) given above are quite natural as they are 
assumptions which do not \emph{a priori} require the existence of an 
invariant measure. However they do imply the existence of an 
\emph{acip} \( \mu \) 
which has all positive Lyapunov exponents. Thus the system 
\( (f, \mu) \) is also nonuniformly expanding in the more abstract 
sense of \emph{Pesin theory}, see \cite{BarPes04}. 
The systematic construction of induced Markov maps in many examples 
and under quite general assumptions, as described above, naturally 
leads to the question of whether such a construction is always 
possible in this abstract setting. 
Since the existence of an induced Markov map implies 
nonuniform expansivity this would essentially give an equivalent 
characterization of nonuniform expansivity. A general result in this 
direction has been given for smooth one-dimensional maps in 
\cite{San03}. It would be interesting to extend this to arbitrary 
dimension. A generalization to nonuniformly hyperbolic 
surface diffeomorphisms is work in progress \cite{LuzSan}.

It seems reasonable to believe that the scope of application of 
induced Markov towers may go well beyond the statistical properties 
of 
a map \( f \). The construction of the induced Markov map 
in \cite{San03} for example is 
primarily motivated by the study if the Hausdorff dimension of 
certain sets. A particularly interesting application would be a 
generalization of the existence (parameter exclusion) argument 
sketched in section \ref{existence}. Even in a very general setting, 
with no information about the map \( f \) except perhaps the 
existence 
of an induced Markov map, it is natural to ask about the possible 
existence of induced Markov maps for small perturbations of \( f \). 
If, moreover, the existence of an induced Markov map were essentially 
equivalent to nonuniform expansivity then this would be a question 
about the persistence of nonuniform expansivity under small 
perturbations. 
\begin{conjecture}
    Suppose that \( f \) is nonuniformly expanding. Then sufficiently 
    small perturbations of \( f \) have positive probability of also 
    being nonuniformly expanding. 
 \end{conjecture}
 Using the Markov induced maps one could define, even in a very 
 abstract setting, a \emph{critical region} \( \Delta \) formed by 
 those points that have very large return time. Then outside \( 
\Delta \) 
 one would have essentially uniform expansivity and these, as well as 
 the Markov structure, would essentially persist under small 
 perturbations.  One could then perturb \( f \) and, up to parameter 
 exclusions, try to show that the Markov structure can be extended 
 once again to the whole of \( \Delta \) for some nearby map \( g \). 

 \subsection{Verifying nonuniform expansivity}
\index{Nonuniform expansivity ! verifying}
The verification of the conditions of nonuniform hyperbolicity are a 
big problem on both a theoretical and a practical level.  
 As mentioned in Section \ref{existence}, for the important class of 
 one-dimensional maps with critical point, nonuniform expansivity 
 occurs with positive probability but for sets of parameters which 
are 
 topologically negligible and thus essentially impossible to pinpoint 
 exactly. The best we can hope for is to show they occur with ``very 
 high'' probability in some given small range of parameter values.
 
 However even this is generally impossible with the available 
 techniques. 
Indeed, all existing argument rely on choosing a \emph{sufficiently 
small} parameter 
 interval \emph{centred} on some \emph{sufficiently good}
 parameter value. The 
 closeness to this parameter value is then used to obtain the various 
 conditions which are required to start the induction.  
 However the problem then reduces to showing that such a good 
 parameter value exists in the particular parameter interval of 
 interest, and this is again both practically and theoretically 
impossible 
 in general. Moreover, even if such a parameter value was determined 
 (as in the special case of the ``top'' quadratic map) existing 
 estimates do not control the size of the neighbourhood in which the 
 good parameters are obtained nor the relative proportion of good 
 parameters. For example there are no explicit bounds for the actual 
 measure of the set of parameters in the quadratic family which have 
 an \emph{acip}. A standard coffee-break joke directed towards 
 authors of the papers on the existence of such maps is that so much 
 work has gone into proving the existence of a set of parameters 
 which as far as we know might be infinitesimal. Moreover, there just 
 does not seem to be any even heuristic argument for believing that 
 such a set is or isn't very small. Thus, for no particular reason 
 other than a reaction to these coffee-break jokers (!), 
 we formulate the  following 
 \begin{conjecture}
     The set of parameters in the quadratic family which admit an acip
     is ``large''.
     \end{conjecture}
  For definiteness let us say that ``large'' means at least \( 50\% 
\) but 
  it seems perfectly reasonable to expect even \( 80\% \)   or \( 
90\% 
  \), and of course we mean here those parameters between the 
  Feigenbaum period doubling limit and the top map. 
  An obvious strategy for proving (or disproving) this conjecture 
would 
be to develop a technique for estimating the proportion of maps 
having 
an \emph{acip} in any given small one-parameter family of maps. The 
large parameter interval of the quadratic family could then be 
subdivided into small intervals each and the contribution of each of 
these small intervals could then be added up. 

A general technique of this kind would also be interesting in a much 
broader 
context of applications. As mentioned in the introduction, 
many real-life systems appear to have a combination of deterministic 
and random-like behaviour which suggests that some form of 
expansivity and/or hyperbolicity might underly the basic driving 
mechanisms. In modelling such a system it seems likely that one may 
obtain a parametrized family and be interested in a possibly narrow 
range of parameter values. It would be desirable therefore to be able 
to obtain a rigorous prove of the existence of stochastic like 
behaviour such as mixing with exponential decay of correlations in 
this family and to be able to estimate the probability of such 
behaviour occurring. An extremely promising strategy has been 
proposed 
recently by K. Mischaikow. The idea is to combine non-trivial 
numerical estimates with the geometric and probabilistic parameter 
exclusion argument discussed above. Indeed the parameter exclusion
argument, see section \ref{existence}, relies 
fundamentally on an induction which shows that the probability of 
being excluded at time \( n \) are exponentially small in \( n \). 
The implementation of this argument however also requires several 
delicate relations between different system constants to be satisfied 
and in particular no exclusions to be required before some 
sufficiently large \( N \) so that the exponentially small exclusions 
occurring for \( n>N \) cannot cumulatively add up to the full 
measure 
of the parameter interval under consideration. The assumption of the 
existence of a 
particularly good parameter value \( a^{*} \) and the assumption that 
the parameter interval is a sufficiently small neighbourhood of \( 
a^{*} 
\) are used in all existing proofs to make sure that certain 
constants can be chosen arbitrarily small or arbitrarily large thus 
guaranteeing that the necessary relations are satisfied. Mischaikow's 
suggestion is to reformulate the induction argument in such a way 
that 
the inductive assumptions can, at least in principle, be explicitly 
verified computationally. This requires the dependence of all the 
constants in the 
argument to be made completely explicit in such a way that the 
inductive assumptions boil down to a finite set of open conditions on 
the family of maps which can be verified with finite precision in 
finite time. 

Besides the interest of the argument in this particular 
setting this could perhaps develop into an extremely fruitful 
interaction between the ``numerical'' and the 
``geometric/probabilistic''
 approach to Dynamical Systems, and contribute significantly to the 
 applicability of the powerful methods of Dynamical Systems to the 
 solution and understanding of real-life phenomena.

\printindex
\end{document}